# An Interior Penalty coupling strategy for Isogeometric non-conformal Kirchhoff-Love shell patches

Giuliano Guarino[1,2], Pablo Antolin[1], Alberto Milazzo[2], and Annalisa Buffa[1]

[1]Institute of Mathematics, École Polytechnique Fédérale de Lausanne, CH-1015, Lausanne, Switzerland

[2]Department of Engineering, Università degli Studi di Palermo, 90128, Italy

**Abstract**

This work focuses on the coupling of trimmed shell patches using Isogeometric Analysis, based on higher continuity splines that seamlessly meet the $C^1$ requirement of Kirchhoff-Love-based discretizations. Weak enforcement of coupling conditions is achieved through the symmetric interior penalty method, where the fluxes are computed using their correct variationally consistent expression that was only recently proposed and is unprecedentedly adopted herein in the context of coupling conditions. The constitutive relationships account for generically laminated materials, although the proposed tests are conducted under the assumption of uniform thickness and lamination sequence. Numerical experiments assess the method for an isotropic and a laminated plate, as well as an isotropic hyperbolic paraboloid shell from the new shell obstacle course. The boundary conditions and domain force are chosen to reproduce manufactured analytical solutions, which are taken as reference to compute rigorous convergence curves in the $L^2$, $H^1$, and $H^2$ norms, that closely approach optimal ones predicted by theory. Additionally, we conduct a final test on a complex structure comprising five intersecting laminated cylindrical shells, whose geometry is directly imported from a STEP file. The results exhibit excellent agreement with those obtained through commercial software, showcasing the method's potential for real-world industrial applications.

## 1 Introduction

Due to their ability to efficiently distribute stress within their volume, shell structures are extensively utilized in high-performance industrial applications. With the advent of composite materials, laminates, which are created by layering orthotropic materials with fibers oriented to optimize the mechanical response, have garnered significant attention from the scientific community. Experimental and numerical investigations of such laminates have been crucial in meeting design criteria and enhancing overall structural performance.

The numerical investigation of laminates using full-scale three-dimensional analysis in practical application is often avoided due to its extensive computational requirements, especially in the initial stages of design. As an alternative, two-dimensional versions of theories describing the main mechanical phenomena involving thin-walled structures have been proposed. The classical Kirchhoff-Love theory [1, 2, 3] assumes that a unit segment perpendicular to the shell's mid-surface remains straight and perpendicular to the surface. The Reissner-Mindlin theory [4, 5], also known as the First-order Shear Deformation Theory (FSDT), relaxes the assumption of perpendicularity and becomes more suitable for moderately thick shells. Additionally, Higher-Order models [6, 7] consider more complex displacement behaviors along the thickness, constituting a sort of intermediate level in terms of accuracy between two-dimensional and three-dimensional theories.

The Kirchhoff-Love shell theory requires only the displacement field in the mid-surface of the shell as main variable, since the the rotation of the perpendicular unit segment can be computed directly from the derivatives of the latter. However, due to the fourth-order nature of the equations, the continuity of the approximation space needs to be $C^1$ over the shell mid-surface, meaning that both the displacements and their first derivatives need to be continuous across elements' boundaries. This requirement poses challenges in the choice of basis functions and the construction of the approximation space, especially in the context of finite element analysis where Lagrange-type basis functions may not satisfy the $C^1$ continuity condition, resulting in a preference for shell elements based on Reissner-Mindlin in the available commercial software. Moreover, when the shell mid-surface is characterized by kinks, to preserve the angle between intersecting faces, the condition on the continuity of the rotation of the perpendicular unit segment along the common edge does not coincide with the $C^1$ continuity of the displacement field anymore, requiring additional effort to be enforced. Early approaches to solve the Kirchhoff-Love equations in the context of Finite Element Method relied on $C^0$ continuity of the basis functions, with $C^1$ continuity enforced weakly in a mixed Continuous-Discontinuous Galerkin approach [8].

The choice of the most suitable 2D theory and the numerical method to solve the resulting set of partial differential equations (PDEs) is an active research topic, with various methods being explored. Examples of these methods include the Finite Element Method [9, 10, 11], Discontinuous Galerkin method [12, 13, 14], Ritz method [15, 16], or Differential Quadrature method [17, 18].

A particularly promising approach that recently proposed is the Isogeometric Analysis (IGA). In IGA, NURBS (Non-Uniform Rational B-Splines) basis functions are used both to define the surface of the shell and to construct the approximation space for the primary variables, allowing for a seamless connection between design and analysis [19]. Since its introduction in [20], IGA has been successfully applied to solve Kirchhoff-Love [21, 22], Reissner-Mindlin [23, 24, 25], and Higher-Order [26] shell theories. In particular, since NURBS functions can be easily constructed with arbitrary continuity, satisfying the $C^1$ requirement for the Kirchhoff-Love shell equation becomes straightforward withing each IGA patch.

However, when dealing with complex shapes, multiple IGA patches are often required to accurately represent the desired geometry, and efficiently coupling these patches becomes a critical issue. Various approaches have been proposed to couple adjacent IGA patches in a strong sense by directly linking some of their degrees of freedom: For example, in [27], IGA regions of the domain are connected in a strong sense with other regions modeled with a mesh-free approach; in [28] patches meeting at $G^0$ interfaces are both connected in a strong sense to auxiliary bending strips that approximate the kink; in [29], the approximation functions for the displacement are continuous across the patches with only the coupling of the rotation imposed in a weak sense; in [30], the coupling approach relies on the Reissner-Mindlin theory where also rotation degrees of freedom are directly available; the construction of $C^1$ multi-patch approximation spaces, as detailed in [31, 32], is used in [33] and [34] for the strong coupling of both displacement and rotation for Kirchhoff-Love IGA patches, but limited to $G^1$ geometries. However, all of these approaches rely on a conforming requirement, meaning that the parameterization of the common edge is the same for each of the patches to be coupled.

For discretization involving IGA patches meeting at non-conforming interfaces, the continuity of the displacement and rotation must be enforced with a weak method. In fact, the main advantage of such approach is that the coupling condition, as well as the boundary Dirichlet condition, does not need to be intrinsically satisfied by the solution space, allowing more flexibility in its definition. In the literature, various methods have been proposed to weakly enforce coupling between IGA patches for many model problems, including the Kirchhoff-Love equations. Examples of these methods include the morthar type and Lagrange multipliers methods [35, 36, 37, 38, 39, 40, 41], the pure penalty methods [42, 43, 44, 45, 46, 47], the projected super-penalty method [48, 49] and the Nitsche-type methods [50, 51, 52, 53, 54, 55, 56, 57, 58]. Among these, the Nitsche-type methods are particularly appealing as they do not require the introduction of additional degrees of freedom as in the Lagrange multipliers methods and, when properly stabilized, do not suffer from the ill-conditioning issues typically seen in penalty approaches. However, constructing a Nitsche-type method for the Kirchhoff-Love shell equations requires computing the fluxes for the formulation. The expression typically found in the literature [59, 51], that tracks back to Koiter's work [3], was recently identified as incorrect in [60]. In that study, a new expression was proposed and validated by several rigorous numerical tests related to the weak enforcement of essential boundary conditions.

Indeed, while multi-patch NURBS offer the potential to construct geometries of any desired curvature profile, the design of complex structures involving multiple intersecting surfaces, cut-outs, or local features can often result in what is commonly referred to as dirty geometries, where boundaries between patches are not watertight. To overcome these limitations, a branch of research in IGA has focused on developing spline spaces using more intricate unstructured grids. The underlying idea is that by enhancing the flexibility of the spline space to accommodate mesh topologies different from the classical one based on a tensor-product rectangular grid, these technologies can serve as tools during the design phase to generate surfaces with more diverse shapes, enabling the representation of local features with high resolution. Numerous instances in this direction have been proposed in recent literature. To cite some, in the context of Kirchhoff-Love shells, the work presented in [61] reparameterizes surfaces with cut-outs using third-degree analysis-suitable T-splines (AST-splines) that are $C^2$ everywhere. In [62], multiple rectangular T-splines are combined through the use of extraordinary points where continuity is locally reduced to $C^1$. AST-splines are further extended in [63] to ensure non-negativity of the basis functions and in [64] to allow for multiple extraordinary points per face. In [65], manifold-based basis functions with selected $C^0$ edges are proposed to model shells with kinks, where the continuity of the rotation is locally enforced through pure penalty. Unstructured splines (U-splines) are introduced in [66] capable of connecting rectangular and triangular cells, retaining higher continuity except in the triangles' edges. Fully triangle configuration B-splines (TCB-splines) are instead adopted in [67] to reparameterize and analyze shell structures topologically equivalent to a disk with a finite number of holes. In [68], the G-spline technology allows for building splines with an arbitrary unstructured quadrilateral layout while maintaining global $C^1$ continuity. In these approaches, to some extent, there is a departure from the traditional concept of designing by assembling different patches. Although some companies in the field are beginning to implement unstructured spline technologies in commercial software, their widespread adoption in industrial applications would necessitate a shift in the design paradigms of practitioners, which might require time and effort. Furthermore, shell mid-surfaces with kinks impose a $C^0$ continuity on the approximation space, still necessitating a weak coupling for the rotation.

As an alternative, the trimmed approach also allows for the definition of complex surfaces while limiting the number of NURBS patches required. This involves embedding a trimming curve in the parametric domain of a surface to delimit its outer boundary and identify active and non-active regions. However, this increased flexibility comes with certain challenges. In fact, when adapting the IGA paradigm to this approach for representing geometries, the presence of trimmed elements raises issues related to integration, conditioning of the linear system, and stability of the method. These same issues arise when coupling between two patches occurs at a trimmed boundary, where one or both of the patches may be trimmed by the interface [54, 48]. Addressing these challenges is crucial to ensure the robustness and accuracy of the method.

Regarding the integration over trimmed elements, several techniques have been proposed in the literature, including: i) hierarchical finite cells: where trimmed elements are subdivided into a hierarchy of smaller cells where a standard integration rule is applied [54, 53]; ii) level-set function: that is applicable to domains where the boundary is implicitly represented by the zero level-set of a reference function [69, 70]; blending functions: that are used to approximate the geometry of the trimmed element, enabling efficient integration [71, 52]. In particular, recently, a robust and efficient algorithm based on higher-order reparameterization of trimmed elements has been proposed in [72], which allows dealing with explicitly defined domains.

In this work, we focus on the linear elastic static analysis of thin-walled structures using the Kirchhoff-Love shell equations and the IGA approach. Laminated shells are considered, and the formulation for isotropic ones is deduced as a special case. When dealing with structures composed of multiple trimmed patches, the proposed coupling strategy relies on the symmetric Nitsche method, also known as interior penalty method in the context of coupling [73]. This work seeks to reassess the interior penalty method for Kirchhoff-Love shells in light of the revised expression for the fluxes [60], which are employed for the first time here to enforce coupling conditions. To handle integration over trimmed elements robustly, we employ the algorithm implemented in [72]. This comprehensive approach allows for efficient and accurate analysis of complex laminated shell structures with non-conforming trimmed interfaces and boundaries.

The paper is structured as follows: In Sec.2, B-spline functions and their extension to trimmed domains are described; Sec.3 introduces the formulation for the Kirchhoff-Love shell equations, along with details regarding the involved differential geometry; Sec.4 recalls the correct expression for the formulation fluxes presented in [60], followed by the Nitsche's formulation for the boundary conditions and the formulation for the coupling conditions; In Sec.5, the efficiency of the method is demonstrated through comparisons with pure penalty methods for an isotropic and a laminated Kirchhoff plate, as well as an isotropic hyperbolic paraboloid shell from the new shell obstacle course for Kirchhoff-Love [60], obtaining for the first time $L^2$, $H^1$, and $H^2$ convergence curves for a generally-curved shell in a coupling test. Additionally, the method is applied to a complex structure consisting of five intersecting laminated cylindrical shells, and the results are compared with those obtained using commercial software; Finally, Sec.6 presents the conclusions of the study.

# 2 The Isogeometric Analysis method

In this section, an overview of B-splines and NURBS functions is given. In the context of IGA methods, they are utilized to construct both surfaces of shells and to discretize their displacement field. In this article, shell surfaces are eventually represented using a trimmed approach. This approach begins with a simple background surface that follows a tensor product structure. The surface is then trimmed by defining its boundary through some additional curves. The details of this trimming operations are provided in the following sections.

## 2.1 The B-splines functions

Univariate B-splines are created based on a polynomial order $p$ and a knot vector, which is a sequence of non-decreasing knot values $\Xi = \{\xi^1, \xi^2, ..., \xi^{n+p+1}\}$. These parameters are used to construct $n$ basis functions $N_i^p(\xi)$, where $\xi$ is the curvilinear coordinate and $i = 1, 2, ...n$. The Cox-de Boor recursion formula is employed to generate these basis functions [74]. A B-spline curve embedded in $\mathbb{R}^3$ can be constructed by multiplying the basis functions by some control points $\boldsymbol{P}_i \in \mathbb{R}^3$ and summing up as

$$\boldsymbol{\mathcal{F}}(\xi) = \sum_{i=1}^{n} N_i^p(\xi)\boldsymbol{P}_i \ . \tag{1}$$

It is important to highlight the following aspects of B-splines:

i) From the knot vector $\Xi$, the vector $\Theta = [\eta^1, \eta^2, \dots, \eta^r]$ is constructed taking only the consecutive, non-repeating values $\eta^i$ in such a way that in the intervals $[\eta^i, \eta^{i+1}]$ the partition of unity property of the spline basis functions is satisfied. It is worth noting that the number of elements of $\Theta$, here denoted by $r$, depends on the specificity of $\Xi$. A B-spline function is therefore defined piece-wise in the intervals $[\eta^i, \eta^{i+1}]$.

ii) Within each interval, a B-spline function is infinitely differentiable ($C^\infty$). However, at the knots, the continuity is at most $C^{p-1}$ and is reduced of one unity for every repetition of the knot value.

iii) The piece-wise nature of a B-spline naturally leads to a mesh-like structure, where each element corresponds to a different interval.

From univariate B-splines basis functions, their bivariate counterparts are constructed using a tensor product approach as

$$B_{ij}(\xi_1, \xi_2) = N_i^p(\xi_1)N_j^p(\xi_2) \ , \tag{2}$$

where it is assumed the same polynomial degree in both curvilinear direction $\xi_1$ and $\xi_2$. Additionally, $N_i^p(\xi_1)$ and $N_j^p(\xi_2)$ are univariate basis functions constructed from the knot vectors $\Xi^1 = \{\xi_1^1, \xi_1^2, ..., \xi_1^{n+p+1}\}$ and $\Xi^2 = \{\xi_2^1, \xi_2^2, ..., \xi_2^{m+p+1}\}$, respectively. Therefore, a B-spline surface is constructed as

$$\boldsymbol{\mathcal{F}}(\xi_1, \xi_2) = \sum_{i=1}^{n}\sum_{j=1}^{m} B_{ij}(\xi_1, \xi_2)\boldsymbol{P}_{ij} \ , \tag{3}$$

where $\boldsymbol{P}_{ij} \in \mathbb{R}^3$ is a generic control point. The knot vectors in a B-spline identify the parametric domain. When open knot vectors are adopted the parametric domain is defined as $\Omega_0 = [\xi_1^1, \xi_1^{n+p+1}] \times [\xi_2^1, \xi_2^{n+p+1}]$, meaning that $(\xi_1, \xi_2) \in \Omega_0$. The bivariate splines inherit the piece-wise definition property from univariate ones through the tensor product structure. Therefore, a rectangular Bezier grid is identified on $\Omega_0$ and the domain of a generic cell is denoted as $Q_0 = [\eta_1^i, \eta_1^{i+1}] \times [\eta_2^j, \eta_2^{j+1}]$, where $[\eta_\alpha^i, \eta_\alpha^{i+1}]$ denotes the $i$-th interval of definition of the univariate basis functions corresponding to $\xi_\alpha$.

For more details on 1D and 2D B-splines, together with the extension to NURBS curves and surfaces, that is not reported here for the sake of conciseness, the interest reader is referred to [74, 75].

## 2.2 Space of trimmed splines

A relatively simple approach to represent complex geometries with intricate boundaries and/or internal holes consists in using a standard B-spline surface and delimiting the actual geometry through some simply-connected curves that define the internal and external boundaries. By sampling a sufficient number of points in the physical space, each of these curves is projected onto the parametric domain, where they delineate two regions. The trimming operation selects only one of these two regions. By repeating this operation for each curve, an active subset of the parametric domain $\Omega \subset \Omega_0$ is identified, which maps through $\boldsymbol{\mathcal{F}}(\xi_1, \xi_2)$ to the final surface $S$ of the modeled shell (see Fig. 1). The space of B-splines in the trimmed domain is defined as:

$$S_h = \text{span}\{B_{ij} \circ \boldsymbol{\mathcal{F}}^{-1} : i \in \{1, ..., n\},\ j \in \{1, ..., m\},\ \text{supp}\{B_{ij}\} \cap \Omega \neq 0\} \tag{4}$$

The trimming operation modifies the rectangular Bezier grid by classifying its cells, denoted as $Q$, into active or non-active, accordingly to whether the corresponding untrimmed cell satisfy $Q_0 \cap \Omega \neq \emptyset$. Active elements are further categorized as entire or partial depending on whether $Q$ is equal to the untrimmed cell $Q_0$ or only a portion of it. Partial cells are cut by the trimming curves and require special treatment for integration. In this work, the algorithm presented in [72] is employed, which performs a reparameterization of the cut cells, enabling the application of Gaussian rules and ensuring high-order accuracy in the integration.

# 3 The Kirchhoff-Love shell theory

Unlike other shell theories, the Kirchhoff-Love equation solely considers the displacement of the mean surface of the shell as the primary variable. This theory is based on the assumption that straight segments perpendicular to the mean surface remain both straight and perpendicular after deformation. Consequently, its rotation is directly obtained from the derivatives of the displacement field, resulting in a fourth-order problem that necessitates $C^1$ continuity of the variable. For a more comprehensive understanding of the derivation of the Kirchhoff-Love shell theory, interested readers are referred to [21, 60]. In this section, for the sake of completeness, we present the theory first for a single patch shell and subsequently extend it to a multi-patch setting.

## 3.1 Problem setting

Let $S \in \mathbb{E}^3$ be the mean oriented surface of the shell under consideration, with $\mathbb{E}^3$ representing the three-dimensional Euclidean space. Let $S$ be the image of the parametric domain $\Omega \in \mathbb{R}^2$, accordingly to the

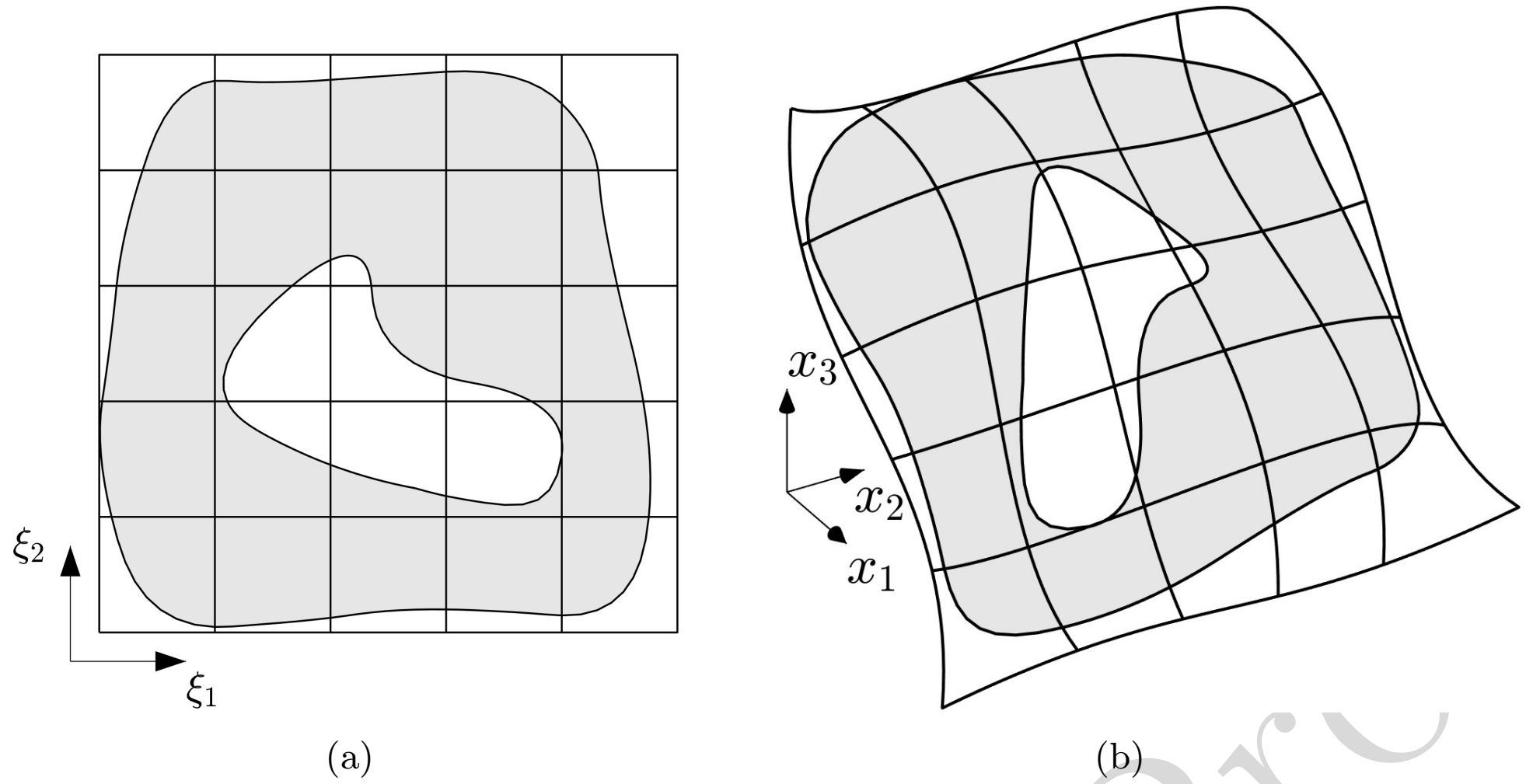

Figure 1: Active (in grey) and non-active (in white) regions of the parametric (a) and physical (b) domains.

map

$$\boldsymbol{x}_0 = \boldsymbol{x}_0(\xi_1, \xi_2) = \begin{bmatrix} x_{01}(\xi_1, \xi_2) \\ x_{02}(\xi_1, \xi_2) \\ x_{03}(\xi_1, \xi_2) \end{bmatrix} , \tag{5}$$

where $\xi_1, \xi_2$ are the curvilinear coordinates spanning $\Omega$. The components of the vector $\boldsymbol{x}_0$ refer to the standard basis $\boldsymbol{e}_1\boldsymbol{e}_2\boldsymbol{e}_3$ of the Euclidean space.

Let $\partial\Omega$ denote the boundary of the parametric domain, which is mapped onto the boundary of the surface $\Gamma$. The boundary $\Gamma$ is split into two parts as $\Gamma = \Gamma^{D_1} \cup \Gamma^{N_1}$, and as $\Gamma = \Gamma^{D_2} \cup \Gamma^{N_2}$ and $\Gamma^{D_1} \cap \Gamma^{N_1} = \emptyset$, $\Gamma^{D_2} \cap \Gamma^{N_2} = \emptyset$. Where $\Gamma^{D_1}$ and $\Gamma^{D_2}$ represent the portions of the boundary where Dirichlet displacement and rotation boundary conditions are applied respectively, while $\Gamma^{N_1}$ and $\Gamma^{N_2}$ represent the portions of the boundary where Neumann force and moment boundary conditions are applied, respectively. Let us define the set of corners $\chi \in \Gamma$ which is further divided into $\chi^D \in \overline{\Gamma^{D_1}}$ where Dirichlet displacement boundary conditions are applied and $\chi^N \in \Gamma^{N_1}$ where Neumann force boundary conditions are applied. The external force distributed on the shell surface is denoted as $\tilde{\boldsymbol{F}}$. Whereas, regarding the surface boundary, the applied force is denoted as $\tilde{\boldsymbol{\tau}}$, while the bending and twisting moments are denoted as $\tilde{M}_{nn}$ and $\tilde{M}_{nt}$, respectively. However, for the Kirchhoff-Love shell theory $\tilde{\boldsymbol{\tau}}$ and $\tilde{M}_{nt}$ cannot be imposed separately and both contribute to the ersatz force applied on $\Gamma^{N_1}$. On the other hand, $\tilde{M}_{nn}$ is the only moment applied on $\Gamma^{N_2}$.

The shell material is assumed to be a laminate, with homogeneous, orthotropic layers having a uniform lamination angle and of uniform thickness. Therefore, the total thickness of the shell, denoted as $\tau$, is also uniform across the surface.

In the reminder of the article, Latin indices span the set $\{1, 2, 3\}$ while Greek indices span the set $\{1, 2\}$. The Einstein summation convention is utilized for repeated indices.

## 3.2 Basics of differential geometry

Let us introduce some concepts of differential geometry needed for the formulation. Starting from the map in Eq.(5), the local covariant basis is defined as

$$\boldsymbol{a}_\alpha = \boldsymbol{x}_{0,\alpha}(\xi_1, \xi_2) \; . \tag{6}$$

Here, the comma preceding one or more Greek indices indicates a series of coordinate derivatives in the specified sequence of curvilinear directions. It is important to note that the vectors of the covariant basis are tangential to the lines of constant curvilinear coordinates and, therefore, lie on the plane that is locally tangent to the surface. As a result, the unit vector $\boldsymbol{a}_3$ locally orthogonal to the surface $S$ is obtained from the covariant basis as

$$\boldsymbol{a}_3 = \frac{\boldsymbol{a}_1 \times \boldsymbol{a}_2}{|\boldsymbol{a}_1 \times \boldsymbol{a}_2|} \ , \tag{7}$$

where $|\bullet|$ is the standard Euclidean norm. The covariant components of the metric tensor are defined as $a_{\alpha\beta} = \boldsymbol{a}_\alpha \cdot \boldsymbol{a}_\beta$, where $\cdot$ denotes the dot product. The determinant of the metric tensor is denoted as $a$. The contravariant components of the metric tensor are obtained from the covariant components as $[a^{\alpha\beta}] = [a_{\alpha\beta}]^{-1}$ and allow us to compute the contravariant basis vectors as

$$\boldsymbol{a}^\alpha = a^{\alpha\beta}\boldsymbol{a}_\beta \ , \tag{8}$$

that satisfy the property $\boldsymbol{a}_\alpha \cdot \boldsymbol{a}^\beta = \delta^\beta_\alpha$, where $\delta^\beta_\alpha$ represents the Kronecker delta. Additionally, we introduce the covariant components and the mixed components of the curvature tensor, defined respectively as

$$b_{\alpha\beta} = \boldsymbol{a}_3 \cdot \boldsymbol{a}_{\alpha,\beta} \ , \tag{9a}$$

$$b^\alpha_\beta = a^{\alpha\gamma} b_{\gamma\beta} \ , \tag{9b}$$

where consistently with the introduced comma notation for the coordinate derivative $\boldsymbol{a}_{\alpha,\beta} = \partial \boldsymbol{a}_\alpha / \partial \xi_\beta = \partial^2 \boldsymbol{x}_0 / \partial \xi_\alpha \partial \xi_\beta$. In this article it is adopted the convention of indicating covariant coordinates referring to the contravariant basis with upper indices, as in $v^\alpha$, and contravariant coordinates referring to the covariant basis with lower indices, as in $v_\alpha$. Furthermore, the notation $v_{\alpha|\beta}$ denotes the covariant derivative of the $\alpha$-th component of a generic vector $\boldsymbol{v}$ along the direction $\beta$. This notation is extended to tensors as well, where $\tau_{\alpha\beta|\gamma}$ represents the covariant derivative of the $\alpha\beta$ component of a generic tensor $\boldsymbol{\tau}$ with respect to the $\gamma$ direction.

### 3.3 Weak form of the Kirchhoff-Love shell equations

Focusing on a single patch shell, the weak form of the Kirchhoff-Love equation is stated as: find $\boldsymbol{u} \in V^u$ such that

$$a(\boldsymbol{u}, \boldsymbol{v}) = f(\boldsymbol{v}) \quad \forall \boldsymbol{v} \in V^v \ , \tag{10}$$

where the choice of the vector spaces $V^v \in H^2(\Omega)$ and $V^u \in H^2(\Omega)$ depends on the specific boundary conditions of the problem. In the discretized version of Eq.(10) the choice of the spaces also takes into account whether the boundary condition are applied in a strong or a weak sense, as explained in Sec.4. The bilinear and the linear forms in Eq.(10) are defined as

$$a(\boldsymbol{u}, \boldsymbol{v}) = \int_S \boldsymbol{\varepsilon}(\boldsymbol{v}) : \boldsymbol{N}(\boldsymbol{u}) \mathrm{d}S + \int_S \boldsymbol{\kappa}(\boldsymbol{v}) : \boldsymbol{M}(\boldsymbol{u}) \mathrm{d}S \ , \tag{11a}$$

$$f(\boldsymbol{v}) = \int_S \boldsymbol{v} \cdot \tilde{\boldsymbol{F}} \mathrm{d}S + \int_{\Gamma^{N_1}} \boldsymbol{v} \cdot \tilde{\boldsymbol{T}} \mathrm{d}\Gamma + \int_{\Gamma^{N_2}} \theta_n(\boldsymbol{v}) \tilde{M}_{nn} \mathrm{d}\Gamma + \sum_{C \in \chi^N} \left( v_3 \tilde{R} \right)\Big|_C \ , \tag{11b}$$

where $\boldsymbol{\varepsilon}$ and $\boldsymbol{\kappa}$ represent the membrane and bending strains, respectively, and their expression is provided in Appendix A, while $\boldsymbol{N}$ and $\boldsymbol{M}$ represent the conjugate generalized force and moment, respectively. These quantities are all rank-2 tensors, and their components are related through the following constitutive equations:

$$N^{\alpha\beta} = \mathbb{A}^{\alpha\beta\gamma\delta} \varepsilon_{\gamma\delta} + \mathbb{B}^{\alpha\beta\gamma\delta} \kappa_{\gamma\delta} \ , \tag{12a}$$

$$M^{\alpha\beta} = \mathbb{C}^{\alpha\beta\gamma\delta} \varepsilon_{\gamma\delta} + \mathbb{D}^{\alpha\beta\gamma\delta} \kappa_{\gamma\delta} \ , \tag{12b}$$

where the coefficients introduced constitute the components of the generalized stiffness tensors for a Kirchhoff-Love laminated shell. Their values depend on both the material and the geometry of the shell. A comprehensive description of how these coefficients are obtained is provided in Appendix B for the sake of completeness. The components of the membrane and bending strains, as well as the normal or bending rotation $\theta_n$, are derived as linear combinations of the first and second coordinate derivatives of the displacement vector. The specific expressions for these quantities can be found in Appendix A.

Regarding the terms related to the applied forces in Eq.(10), in addition to the bending moment $\tilde{M}_{nn}$ and the surface force $\tilde{\boldsymbol{F}}$ that were already introduced in Sec.3.1, two additional forces are introduced: the ersatz forces $\tilde{\boldsymbol{T}}$ and the corner forces $\tilde{R}$. These are defined, respectively, as

$$\tilde{\boldsymbol{T}} = \left(\tilde{\tau}_\alpha - \tilde{M}_{nt} b_{\alpha\beta} t^\beta\right) \boldsymbol{a}^\alpha + \left(\tilde{\tau}_3 + \frac{\partial \tilde{M}_{nt}}{\partial t}\right) \boldsymbol{a}^3 \ , \tag{13}$$

$$\tilde{R} = \lim_{\epsilon \to 0} \left(\tilde{M}_{nt}(\boldsymbol{x} + \epsilon \boldsymbol{t}) - \tilde{M}_{nt}(\boldsymbol{x} - \epsilon \boldsymbol{t})\right) \ , \tag{14}$$

where $t^\alpha$ represents the $\alpha$-th contravariant coordinate of the vector $\boldsymbol{t}$, which is the unit vector locally tangent to the counterclockwise-oriented boundary $\Gamma$. Additionally, $\boldsymbol{n} = \boldsymbol{t} \times \boldsymbol{a}_3$ is the outer unit vector orthogonal to the boundary and lying in the local plane tangent to $S$. The components of the applied moment $\tilde{M}_{nn}$ and $\tilde{M}_{nt}$ are referred to the basis formed by the vectors $\boldsymbol{n}$ and $\boldsymbol{t}$.

At this point it is worth mentioning that Kirchhoff-Love shells, unlike Reissner-Mindlin ones, do not suffer from shear locking. However, membrane locking can be an issue in certain critical situations. Various remedies for addressing membrane locking have been proposed in the literature. To mention a few: in [76], two approaches are proposed: the discrete strain gap method, that is based on ad-hoc integration rules for strains, but due to the higher-order continuity of spline spaces in IGA, leads to a loss of sparsity of the stiffness matrix; and a mixed-method incorporating both displacement and membrane stress as main variables in the Hellinger-Reissner principle. The mixed displacement method is introduced in [77], where additional displacement degrees of freedom are introduced from the same approximation space as the actual displacements to compute the assumed membrane strain. Instead, in [78], assumed membrane strains are obtained through a local projection of the B-spline space and then reconstructed on a patch level using a local smoothing procedure. In [79], $C^0$ continuous assumed strains are constructed as a bilinear interpolation of the strains at the four nodes of a rectangular element, addressing locking without the introduction of additional degrees of freedom, the need for additional matrix inversion, and preserving the sparsity of the stiffness matrix, although limited to quadratic elements. It is important to mention that in the present formulation, no treatment for locking is utilized, and an efficient combination of locking treatments with the coupling strategy proposed in this paper is left for further development of the present work.

## 3.4 Extension to multi-patch

The variational statement in Eq.(10) applies to shell structures consisting on a single patch. However, if the structure is composed of $N_P$ patches that intersect at $N_I$ interfaces, the problem becomes finding $\boldsymbol{u} \in \boldsymbol{\mathcal{V}}^u$ such that:

$$\sum_{p=1}^{N_P} a^p(\boldsymbol{u}, \boldsymbol{v}) = \sum_{p=1}^{N_P} f^p(\boldsymbol{v}) \quad \forall \boldsymbol{v} \in \boldsymbol{\mathcal{V}}^v \ , \tag{15}$$

where $a^p(\boldsymbol{u}, \boldsymbol{v})$ and $f^p(\boldsymbol{v})$ are defined as in Eq.(11), but a superscript $p$ is added to indicate that they belong to the $p$-th patch. Additionally, $\boldsymbol{\mathcal{V}}^u$ and $\boldsymbol{\mathcal{V}}^v$ are vector spaces defined over the union of the surfaces of the patches of the structure. These spaces have to be defined in a way that, apart from satisfying the essential boundary conditions, ensures the following conditions over each of the interfaces $\Gamma^i$:

$$[\boldsymbol{u}] = \boldsymbol{0} \ , \tag{16a}$$

$$[\theta_n] = 0 \ , \tag{16b}$$

where $[\bullet]$ represents the jump operator, which calculates the difference between the quantity of interest computed from the different patches at the interface.

In order to enforce this coupling condition in a strong sense, Eq.(16) should be embedded in the spaces $\boldsymbol{\mathcal{V}}^v$ and $\boldsymbol{\mathcal{V}}^u$. In practice, after discretization, it is relatively easy to enforce Eq.(16a) in the case where the patches are conforming at the interface, but is not straightforward to do the same for Eq.(16b) and this is restricted to $G^1$ surfaces [31, 33, 80]. However, in situations where the interface is generated in a non-conforming manner, such as when two patches meet on a trimmed boundary, the coupling conditions can only be enforced in a weak sense. The same applies to essential boundary conditions on the trimmed boundary of each patch.

# 4 The interior penalty coupling for Kirchhoff-Love shells

In this section, the problem stated in Eq.(10) in a continuous framework is discretized by selecting appropriate spaces $V^u$ and $V^v$ for a single patch problem and $\boldsymbol{\mathcal{V}}^u$ and $\boldsymbol{\mathcal{V}}^v$ for a multi-patch problem. In an IGA approach, these spaces are constructed starting from the trimmed B-spline space introduced in Sec.2.2. In this article, Dirichlet boundary conditions are applied in either a weak or a strong sense, depending on the specific problem being investigated.

For a single patch shell, in order to impose Dirichlet boundary conditions in a strong sense, it is necessary to ensure that the test and trial functions $\boldsymbol{u}_h$ and $\boldsymbol{v}_h$ of $V_h^u$ and $V_h^v$, respectively, satisfy $\boldsymbol{u}_h = \tilde{\boldsymbol{u}}$ and $\boldsymbol{v}_h = \boldsymbol{0}$ on $\Gamma_h^{D_1}$ and $\theta_n(\boldsymbol{u}_h) = \tilde{\theta}_n$ and $\theta(\boldsymbol{v}_h) = 0$ on $\Gamma_h^{D_2}$. These conditions can be imposed strongly only on conforming edges. However, if Dirichlet boundary conditions are enforced in a weak sense, these requirements no longer need to be satisfied. The same principles apply when considering shell structures composed of multiple patches. More specifically, regarding the coupling conditions, this work adopts exclusively a weak imposition, as explained in the following sections. Additionally, similarly to other fourth-order equations, in the Kirchhoff-Love one there is a continuity requirement on the trial functions of at least $C^1$ [60, 81].

Different methods are available in literature to apply boundary conditions and coupling conditions in a weak sense. In this work, the symmetric Nitsche method [60], also known as the interior penalty method [73] in the context of coupling conditions, is employed. This method requires the computation of the fluxes of the formulation obtained as explained in Sec.4.1.

## 4.1 The fluxes for the Kirchhoff-Love problem

Due to the complexity of the Kirchhoff-Love shell equations, the computation of the fluxes is not a trivial task. In fact, the initially proposed expression for the fluxes by Koiter [3] was found to be incorrect, as discussed in [60]. For a complete derivation of their expression, interested readers are referred to this source. However, the correct definition is also provided here for the sake of completeness.

Recalling the outer-facing unit normal $\boldsymbol{n}$ introduced in Sec.3.3, in this direction, the fluxes associated with the problem are of two types: those corresponding to the ersatz force and those corresponding to the bending moment. Respectively, their definitions are

$$\boldsymbol{T}(\boldsymbol{u}) = T^\alpha \boldsymbol{a}_\alpha + T^3 \boldsymbol{a}_3 \ , \tag{17a}$$

$$M_{nn}(\boldsymbol{u}) = M^{\alpha\beta} n_\alpha n_\beta \ , \tag{17b}$$

where $n_\alpha$ is the component of $\boldsymbol{n}$ referred to $\boldsymbol{a}^\alpha$, and the components of the vector $\boldsymbol{T}$ are defined as

$$T^\alpha = N^{\alpha\beta} n_\beta - b^\alpha_\gamma M^{\gamma\beta} n_\beta - M_{nt} b^\alpha_\gamma t^\gamma \ , \tag{18a}$$

$$T^3 = M^{\alpha\beta}_{|\beta} n_\alpha + (M^{\alpha\beta} n_\alpha t_\beta)_{,t} \ , \tag{18b}$$

where with the notation $(\bullet)_{,t}$ we denote the arc-lenght derivative along the curve that identifies $\boldsymbol{t}$. Additional details on the computation of the terms appearing in Eq.(18) are given in Appendix C.

## 4.2 The Nitsche's method for weakly imposing Dirichlet boundary conditions

In this article, when a strong imposition of Dirichlet boundary conditions is not possible, the symmetric Nitsche's method [60] is employed. For a single-patch shell, the discretized Kirchhoff-Love variational statement is formulated as follows: find $\boldsymbol{u}_h \in V_h^u$ such that

$$a_h(\boldsymbol{u}_h, \boldsymbol{v}_h) + a_n(\boldsymbol{u}_h, \boldsymbol{v}_h) + a_s(\boldsymbol{u}_h, \boldsymbol{v}_h) = f_h(\boldsymbol{v}_h) + f_n(\boldsymbol{v}_h) + f_s(\boldsymbol{v}_h) \quad \forall \boldsymbol{v}_h \in V_h^v \,. \tag{19}$$

Where $V_h^u$, and $V_h^v$ are constructed from the space of B-splines defined over $\Omega$ and their exact definition depends on whether the boundary conditions are imposed strongly in a portion of the boundary. In the equations presented in this section, it is assumed that essential boundary conditions are applied in a weak sense in the entire $\Gamma_h^{D_1}$ and $\Gamma_h^{D_2}$, whereas in the results presented in Sec.5 it is preferred, when possible, to enforce them in a strong sense. Both the bilinear form on the left-hand side and the linear form on the right-hand side of Eq.(19) are constructed by summing three contributions. As for the terms with a subscript $h$, they are defined as follows:

$$a_h(\boldsymbol{u}_h, \boldsymbol{v}_h) = \int_{S_h} \boldsymbol{\varepsilon}(\boldsymbol{v}_h) : \boldsymbol{N}(\boldsymbol{u}_h) \mathrm{d}S + \int_{S_h} \boldsymbol{\kappa}(\boldsymbol{v}_h) : \boldsymbol{M}(\boldsymbol{u}_h) \mathrm{d}S \,, \tag{20a}$$

$$f_h(\boldsymbol{v}_h) = \int_{S_h} \boldsymbol{v} \cdot \tilde{\boldsymbol{F}} \mathrm{d}S + \int_{\Gamma_h^{N_1}} \boldsymbol{v}_h \cdot \tilde{\boldsymbol{T}} \mathrm{d}\Gamma + \int_{\Gamma_h^{N_2}} \theta_n(\boldsymbol{v}_h) \tilde{M}_{nn} \mathrm{d}\Gamma + \sum_{C \in \chi_h^N} \left( v_{3h} \tilde{R} \right)\Big|_C \,, \tag{20b}$$

and constitute the discretized version of Eq.(10). As such, $S_h$, $\Gamma_h^{N_1}$, $\Gamma_h^{N_2}$, and $\chi_h^N$ are the approximated versions after discretization of $\Gamma^{N_1}$, $\Gamma^{N_2}$, and $\chi^N$, respectively. The terms with a subscript $n$ are the symmetric Nitsche terms, which include the fluxes described in Sec.4.1. These terms are defined as

$$\begin{aligned} a_n(\boldsymbol{v}_h, \boldsymbol{u}_h) = &- \int_{\Gamma_h^{D_1}} \left( \boldsymbol{T}(\boldsymbol{v}_h) \cdot \boldsymbol{u}_h + \boldsymbol{v}_h \cdot \boldsymbol{T}(\boldsymbol{u}_h) \right) \mathrm{d}\Gamma - \int_{\Gamma_h^{D_2}} \left( M_{nn}(\boldsymbol{v}_h) \theta_n(\boldsymbol{u}_h) + \theta_n(\boldsymbol{v}_h) M_{nn}(\boldsymbol{u}_h) \right) \mathrm{d}\Gamma + \\ &- \sum_{C \in \chi_h^D} \left( R(\boldsymbol{v}_h) u_{3h} + v_{3h} R(\boldsymbol{u}_h) \right)\big|_C \,, \end{aligned} \tag{21a}$$

$$f_n(\boldsymbol{v}_h) = - \int_{\Gamma_h^{D_1}} \boldsymbol{T}(\boldsymbol{v}_h) \cdot \tilde{\boldsymbol{u}} \mathrm{d}\Gamma - \int_{\Gamma_h^{D_2}} M_{nn}(\boldsymbol{v}_h) \tilde{\theta}_n \mathrm{d}\Gamma - \sum_{C \in \chi_h^D} \left( R(\boldsymbol{v}_h) \tilde{u}_3 \right)\big|_C \,, \tag{21b}$$

where $\tilde{\boldsymbol{u}}$ and $\tilde{\theta}_n$ are the applied displacement and normal rotation, respectively and $\Gamma_h^{D_1}$, $\Gamma_h^{D_2}$, and $\chi_h^D$ are the approximated versions of $\Gamma^{D_1}$, $\Gamma^{D_2}$, and $\chi^D$, respectively. The subscript 3 in $u_{3h}$ and $v_{3h}$ denotes the component of the respective vector relative to $\boldsymbol{a}^3 = \boldsymbol{a}_3$. Additionally, the following definition has been employed:

$$R = \lim_{\epsilon \to 0} \left( M_{nt}(\boldsymbol{x} + \epsilon \boldsymbol{t}) - M_{nt}(\boldsymbol{x} - \epsilon \boldsymbol{t}) \right) \,, \tag{22}$$

where $M_{nt} = M^{\alpha\beta} n_\alpha t_\beta$ and $t_\beta$ is the component of $\boldsymbol{t}$ along $\boldsymbol{a}^\beta$. Both $a_n(\boldsymbol{v}_h, \boldsymbol{u}_h)$ and $f_n(\boldsymbol{v}_h)$ are composed of three terms. The first and second terms correspond to the displacements and rotation boundary conditions, respectively. The third term, which involves the displacements at the corners, is introduced to ensure optimal convergence, as discussed in [60]. Each of the three contributes are constructed from a consistency term (e.g., $\boldsymbol{v}_h \cdot \boldsymbol{T}(\boldsymbol{u}_h)$) and a symmetry term (e.g., $\boldsymbol{T}(\boldsymbol{v}_h) \cdot \boldsymbol{u}_h$ and $\boldsymbol{T}(\boldsymbol{v}_h) \cdot \tilde{\boldsymbol{u}}$). Finally, the stabilization terms in Eq.(19), denoted by a subscript $s$, are defined as

$$a_s(\boldsymbol{v}_h, \boldsymbol{u}_h) = \int_{\Gamma_h^{D_1}} \mu_D^b \boldsymbol{v}_h \cdot \boldsymbol{u}_h \mathrm{d}\Gamma + \int_{\Gamma_h^{D_2}} \mu_R^b \theta_n(\boldsymbol{v}_h) \theta_n(\boldsymbol{u}_h) \mathrm{d}\Gamma + \int_{\Gamma_h^{D_1}} \mu_3^b v_{3h} u_{3h} \mathrm{d}\Gamma + \sum_{C \in \chi_h^D} \left( \mu_C^b v_{3h} u_{3h} \right)\big|_C \,, \tag{23a}$$

$$f_s(\boldsymbol{v}_h) = \int_{\Gamma_h^{D_1}} \mu_D^b \boldsymbol{v}_h \cdot \tilde{\boldsymbol{u}} \mathrm{d}\Gamma + \int_{\Gamma_h^{D_2}} \mu_R^b \theta_n(\boldsymbol{v}_h) \tilde{\theta}_n \mathrm{d}\Gamma + \int_{\Gamma_h^{D_1}} \mu_3^b v_{3h} \tilde{u}_3 \mathrm{d}\Gamma + \sum_{C \in \chi_h^D} \left( \mu_C^b v_{3h} \tilde{u}_3 \right)\big|_C \,, \tag{23b}$$

where $\mu_D^b$, $\mu_R^b$, $\mu_3^b$, and $\mu_C^b$ are the so-called penalty parameters, which play a crucial role in the weak imposition of boundary and coupling conditions, since they provide stability to the method. The choice of these parameters is still an open question and depends on the specific problem at hand. A discussion on the importance of these parameters, some guidelines on how to choose them, and the approach adopted in this paper is presented in Sec.4.4.

### 4.3 The interior penalty method for coupling IGA patches

The framework presented in Sec.4.2 focuses on a single-patch shell. However, when dealing with structures composed of multiple shells that intersect at common interfaces, the formulation needs to be extended to address the coupling conditions. To achieve this, the terms in Eq.(19) are enriched with a superscript $p$ to indicate that they belong to the $p$-th patch. Then, the following forms are defined:

$$a^p(\boldsymbol{u}_h, \boldsymbol{v}_h) = a_h^p(\boldsymbol{u}_h, \boldsymbol{v}_h) + a_n^p(\boldsymbol{u}_h, \boldsymbol{v}_h) + a_s^p(\boldsymbol{u}_h, \boldsymbol{v}_h) \ , \tag{24a}$$

$$f^p(\boldsymbol{v}_h) = f_h^p(\boldsymbol{v}_h) + f_n^p(\boldsymbol{v}_h) + f_s^p(\boldsymbol{v}_h) \ . \tag{24b}$$

The discretized version of the Kirchhoff-Love shell equation for multi-patch structures becomes: find $\boldsymbol{u}_h$ in $\boldsymbol{\mathcal{V}}_h^u$ such that

$$\sum_{p=1}^{N_P} a^p(\boldsymbol{u}_h, \boldsymbol{v}_h) + \sum_{i=1}^{N_I} b^i(\boldsymbol{u}_h, \boldsymbol{v}_h) = \sum_{p=1}^{N_P} f^p(\boldsymbol{v}_h) \quad \forall \boldsymbol{v}_h \in \boldsymbol{\mathcal{V}}_h^v \ , \tag{25}$$

where $\boldsymbol{\mathcal{V}}_h^u$ and $\boldsymbol{\mathcal{V}}_h^v$ are the discretized spaces correspondent to $\boldsymbol{\mathcal{V}}^u$ and $\boldsymbol{\mathcal{V}}^v$, respectively, and $b^i(\boldsymbol{u}_h, \boldsymbol{v}_h)$ is the contribute to the variational statement ensuring the coupling between the patches intersecting at the $i$-th interface. These terms are obtained as

$$b^i(\boldsymbol{u}_h, \boldsymbol{v}_h) = b_n^i(\boldsymbol{u}_h, \boldsymbol{v}_h) + b_s^i(\boldsymbol{u}_h, \boldsymbol{v}_h) \ . \tag{26}$$

Once again, Nitsche and stabilization terms have been introduced, denoted by the subscripts $n$ and $s$, respectively. The definitions of these terms are

$$b_n^i(\boldsymbol{u}_h, \boldsymbol{v}_h) = -\int_{\Gamma_h^i} \left( [\boldsymbol{v}_h] \cdot \{\boldsymbol{T}(\boldsymbol{u}_h)\} + [\theta_n(\boldsymbol{v}_h)]\{M_{nn}(\boldsymbol{u}_h)\} \right) \mathrm{d}\Gamma - \gamma_1 \int_{\Gamma_h^i} \left( \{\boldsymbol{T}(\boldsymbol{v}_h)\} \cdot [\boldsymbol{u}_h] + \{M_{nn}(\boldsymbol{v}_h)\}[\theta_n(\boldsymbol{u}_h)] \right) \mathrm{d}\Gamma \ , \tag{27a}$$

$$b_s^i(\boldsymbol{u}_h, \boldsymbol{v}_h) = \int_{\Gamma_h^i} \left( \mu_D^c [\boldsymbol{v}_h] \cdot [\boldsymbol{u}_h] \right) \mathrm{d}\Gamma + \int_{\Gamma_h^i} \left( \mu_R^c [\theta_n(\boldsymbol{v}_h)][\theta_n(\boldsymbol{u}_h)] \right) \mathrm{d}\Gamma \ . \tag{27b}$$

Where $\Gamma_h^i$ denotes the approximation of the $i$-th interface $\Gamma^i$. The penalty terms $\mu_D^c$ and $\mu_R^c$, associated to the displacement and the rotation coupling conditions, respectively, have been introduced. Their choice is discussed together with their counterparts for boundary conditions in Sec.4.4. $\{\bullet\}$ and $[\bullet]$ are the average and jump operators, defined as

$$\{\bullet\} = \gamma_2 \ \bullet^+ + (1 - \gamma_2) \bullet^- \ , \tag{28a}$$

$$[\bullet] = \bullet^+ - \bullet^- \ , \tag{28b}$$

where $\bullet$ denotes a generic quantity defined over both patches at the same point on the interface $\Gamma^i$. The superscript $+$ and $-$ are used to distinguish between the two patches. The parameters $\gamma_1$ and $\gamma_2$ are used to differentiate amongst the Nitsche type methods. Their meaning and effect on the stabilization of the method is discussed in Sec.4.4. In this contribute, these are chosen as $\gamma_1 = 1$ and $\gamma_2 = 0.5$ that leads to a symmetric interior penalty formulation.

In this formulation, it is assumed that the curves in the Euclidean space describing the $i$-th interface are known explicitly. Consequently, the unit vector $\boldsymbol{t}$ tangent to the interface is also assumed to be known. $\boldsymbol{t}$ denotes the tangent vector to both the external boundaries of the patch (as mentioned in Sec.3.3) and the interface between two patches. Moreover, since an interface is common to both the intersecting patches that generate it, the unit vector $\boldsymbol{t}$ is also the same for both patches. The context always makes it clear which condition is being referred to, whether a boundary or a coupling one. However, in the case of the coupling condition, the requirement of $\boldsymbol{t}$ to be oriented in a counter-clockwise direction is discarded, and whether the orientation is clockwise or counter-clockwise depends on the relative position between $\boldsymbol{t}$, $\boldsymbol{a}_3$, and $\boldsymbol{n}$ for each specific patch, being $\boldsymbol{n}$ the outer unit vector normal to the interface, lying on the plane locally tangent to the patch surface.

As a result, $\boldsymbol{n}^+$ and $\boldsymbol{n}^-$ can be different. In the simplest case, they lie in the same direction, either coinciding or being opposite. But, if the patches meet at an angle, their directions differ. The formulation presented in this article is capable of handling every possible case. However, in order to properly compute the average, the fluxes for the second patch (in contrast to the first patch) are obtained with respect to a normal vector entering the surface domain.

## 4.4 Choice of the parameters of the methods

In order to weakly enforce essential boundary and coupling conditions, various methods have been investigated in literature. In particular, different penalty and Nitsche's methods can be constructed based on the presence of Nitsche and/or stabilization terms, the presence and the sign of the symmetry terms, and the definition of the average operator. The unified formulation proposed in [82] introduces the parameters $\gamma_1$ and $\gamma_2$, which allow for the construction of different Nitsche type methods. When no Nitsche terms are present in the formulation, the resulting method is the pure penalty [83, 42, 43, 44, 45, 47]. The symmetric interior penalty method [51, 52, 53, 57] is constructed adopting $\gamma_1 = 1$ and $\gamma_2 = 0.5$, requiring the introduction of penalty terms to ensure stability. The skew Nitsche's method [54] is obtained with $\gamma_1 = -1$ and $\gamma_2 = 0.5$ without adding stabilization terms. This choice leads to a skew-symmetric solving linear system and has the advantage of being parameter-free. Taking $\gamma_1 = -1$ and $\gamma_2 \neq 0.5$ leads to the weighted non-symmetric Nitsche's method [84, 85, 55, 56]. This method provides increased stability compared to other Nitsche's methods but at the expense of losing any symmetry of the linear system.

In pure penalty methods, which rely solely on penalty integrals, the choice of the penalty parameter can significantly impact the accuracy of the method, setting a lower bound on the achievable error corresponding to that particular value of the penalty parameter [47]. If the value is too low, it results in a weak enforcement of the boundary/coupling condition. Conversely, if the value is too high, it leads to ill-conditioning of the linear system. When employing a Nitsche's method that includes the flux terms, the penalty terms are used solely for stabilization purposes. In this case, the minimum value of the penalty parameter required to achieve optimal convergence is lower than that needed for pure penalty methods, allowing for more flexibility in its selection.

In this work, symmetric Nitsche terms defined in Eq.(21) and Eq.(27a) offer the advantage, compared to other Nitsche's methods, of resulting in symmetric linear systems, which can be beneficial in terms of computational efficiency. However, the interior penalty formulation may lose coercivity, and therefore stability, when using a fixed penalty parameter for discretization including small cut elements. Unfortunately, as of the authors' knowledge, a method stable under all circumstances for coupling trimmed shell patches does not currently exist in the literature. A possible local estimate for the penalty parameter to retain coercitivity relies on solving a local eigenvalue problem [86, 50, 87, 54]. However, this approach can significantly increase computational time and lead to high penalty parameters in some critical scenarios such as coupling of patches with different constitutive properties or coupling of elements with drastically different sizes. In these cases, choosing different values of $\gamma_1$ or $\gamma_2$ can benefit stability while limiting the penalty value, as shown in [88], with the drawbacks of loosing symmetry of the linear system. Fortunately, for the symmetric Nitsche's method, penalty parameters that ensure optimal convergence while preserving a reasonable condition number of the system matrix still span a wide range. The approach chosen here takes advantage of this property and follows the recommendations already available in the literature for the choice of the penalty parameters. Further investigation is needed in this direction, but the development of an unconditionally-stable Nitsche-type coupling method falls beyond the scope of this work.

To properly scale the penalty terms with respect to the problem parameters, a typical construction involves multiplying a problem-independent constant, a problem-dependent term, and a mesh-size-dependent term, that might be raised to a mesh-degree-dependent power, accordingly to the problem at hand and the method adopted. Following [60], and extending the construction to interface coupling, the penalty parameters are chosen as:

$$\mu_D^b = \mu_D^c = \beta E_l \tau / h \ , \tag{29a}$$

$$\mu_R^b = \mu_R^c = \beta E_l \tau^3 / h \ , \tag{29b}$$

$$\mu_3^b = \beta E_l \tau^3 / h^3 \ , \tag{29c}$$

$$\mu_C^b = \beta E_l \tau^2 / h^2 \ , \tag{29d}$$

where $\mu_D^b$ and $\mu_D^c$ are employed for displacement boundary and coupling conditions, respectively. Similarly, $\mu_R^b$ and $\mu_R^c$ are employed for rotation boundary and coupling conditions, respectively. The parameters $\mu_3^b$ and $\mu_C^b$ relate the component of the displacement vector along $\boldsymbol{a}_3$ at the boundaries and at the corners, respectively. The corresponding integrals to these last two penalty terms in Eq.(23) are proven to be necessary for achieving optimal convergence in [60]. In Eq.(29), $E_l$ represents the maximum Young modulus of the laminate (as described in Appendix B), and $h$ is a measure of the mesh size. The problem-independent parameter $\beta$ is taken here as either 10, $10^2$, or $10^3$ for all the penalty terms, as specified in the tests of Sec.5. This choice helps balance the enforcement of the boundary/coupling conditions while maintaining a well-conditioned system.

# 5 Results

In this section, we evaluate the performance of the proposed method through various numerical experiments on benchmark problems involving isotropic and laminated plates and shells. The structures in the proposed tests are loaded by domain forces and are subjected to various boundary conditions, including homogeneous and non-homogeneous Dirichlet conditions. The application of these boundary conditions varies depending on the cases, with some being enforced in a strong manner, while others are applied weakly. When available, the numerical solution are compared with analytical ones. These are manufactured computing the applied domain force from the desired distribution of the displacement field by using the strong form of the Kirchhoff-Love shell equations, as presented in [60], or if the test involves a planar geometry using the Kirchhoff plate equations.

To what concerns software, the open-source MATLAB® library GeoPDEs [89, 90] is utilized, with additional functions implemented for the coupling. High-order integration over trimmed elements and their boundaries is achieved using the algorithm presented in [72], which is based on a reparameterization of the trimmed elements. In particular, integration over non-conforming interfaces requires specific attention. The interfaces must be subdivided in such a way that each curved segment corresponds to a unique element on each of the interface's patches. This segmentation of the interface is then projected onto both parametric domains to compute the quantities of interest from each patch. For a detailed description of this procedure, interested readers are referred to [88].

Notoriously, ill-conditioning is a common challenge in linear systems arising from shell element structures. This issue is further exacerbated by high degree polynomials and weak coupling conditions. Additionally, in immersed boundary approaches, a huge difference in element size can appear in certain refinement level of the discretization, contributing to increase the condition number. To mitigate this issue, in this work it is used a Jacobi preconditioner [91] that, despite its simplicity, has demonstrated remarkable efficacy. Nontheless, developing a more robust preconditioner could be beneficial in addressing this issue in more general cases. However, this potential solution falls outside the scope of the current work and could be pursued as an interesting direction for future research.

The coupling strategy is tested on multi-patch geometries connected at non-conforming trimmed interfaces. The convergence curves for multi-patch discretizations in the $L^2$ norm, $H^1$ seminorm, and $H^2$ seminorm are compared to those correspondent to reference single patch ones. Furthermore, results involving a structure comprised of multiple intersecting laminated cylindrical shells showcase the method's potential for industrial problems.

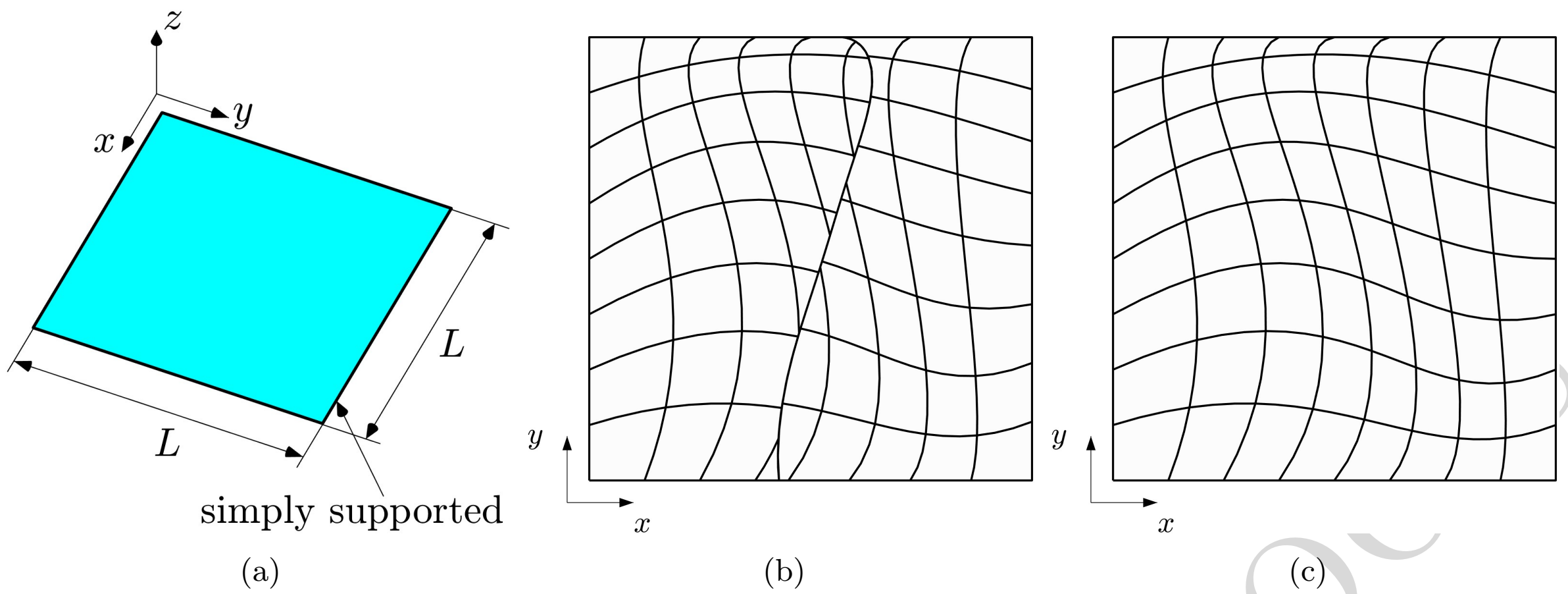

Figure 2: Geometry of the Kirchhoff plate described in Sec.5.1 (a). Discretization of the plate employing two non-conformal IGA patches (b) and a single IGA patch (c).

## 5.1 Square Kirchhoff plate

In this test, the mechanical response of an isotropic plate structure is modelled using two trimmed planar patches coupled along a non-conforming interface. The main geometrical features of the structure are depicted in Fig.2a and the two patches for a specific refinement level are illustrated in Fig.2b. The map of the shell is constructed to ensure that the lines of constant curvilinear coordinates are curved in the physical space, making in such way this example more significant.

To create a non-conforming interface, an additional knot is inserted in each patch, specifically $(0.5, 0.5)$ for the first patch and $(0.45, 0.53)$ for the second patch. After trimming both patches using the same trimming curve, which is also constructed to be curved in the physical domain, they are joint together. This configuration ensures that subsequent dyadic refinements of the discretization maintain the non-conforming nature of the coupling interface.

As a reference for efficiency comparison, Fig.2c displays the same structure modeled with a single patch, which corresponds to the untrimmed left patch of the multi-patch configuration. The material used for the analysis is characterized by a Young modulus $E = 70$ [GPa] and a Poisson ratio $\nu = 0.3$. The plate has a square mid-surface with an edge length $L = 1$ [m] and three values of the thickness are considered, namely $\tau = 0.1$ [m], $\tau = 0.01$ [m], and $\tau = 0.001$ [m]. Simply supported boundary conditions are applied, with homogeneous displacements Dirichlet boundary conditions enforced in a strong manner along the entire external boundaries. A distributed surface force is applied in the direction $\boldsymbol{e}_3$ to reproduce the manufactured smooth solution:

$$\boldsymbol{u}^{ex} = U_0 \sin\left(\frac{2\pi x}{L}\right) \sin\left(\frac{2\pi y}{L}\right) \boldsymbol{e}_3 \, , \tag{30}$$

where $U_0 = 0.1$ [m] is the maximum absolute value of the displacement. Regarding the choice of the penalty parameters, the arbitrary coefficient in Eq.(29) is selected as $\beta = 10^2$.

In Fig.3, the convergence behavior of the $L^2$ norm, $H^1$, and $H^2$ seminorms of the error for different polynomial values, plate thickness, and discretization approaches in the presented test case is illustrated. The dashed lines represent the convergence curves for the single-patch case, serving as a reference for optimal convergence relative to that specific polynomial order. The solid lines correspond to the interior penalty method discussed in Sec.4.3, showing how accurately they follow the reference convergence curves. On the other hand, the dotted lines depict the convergence results for the pure penalty method, where only the stabilization terms in Eq.(27) are considered.

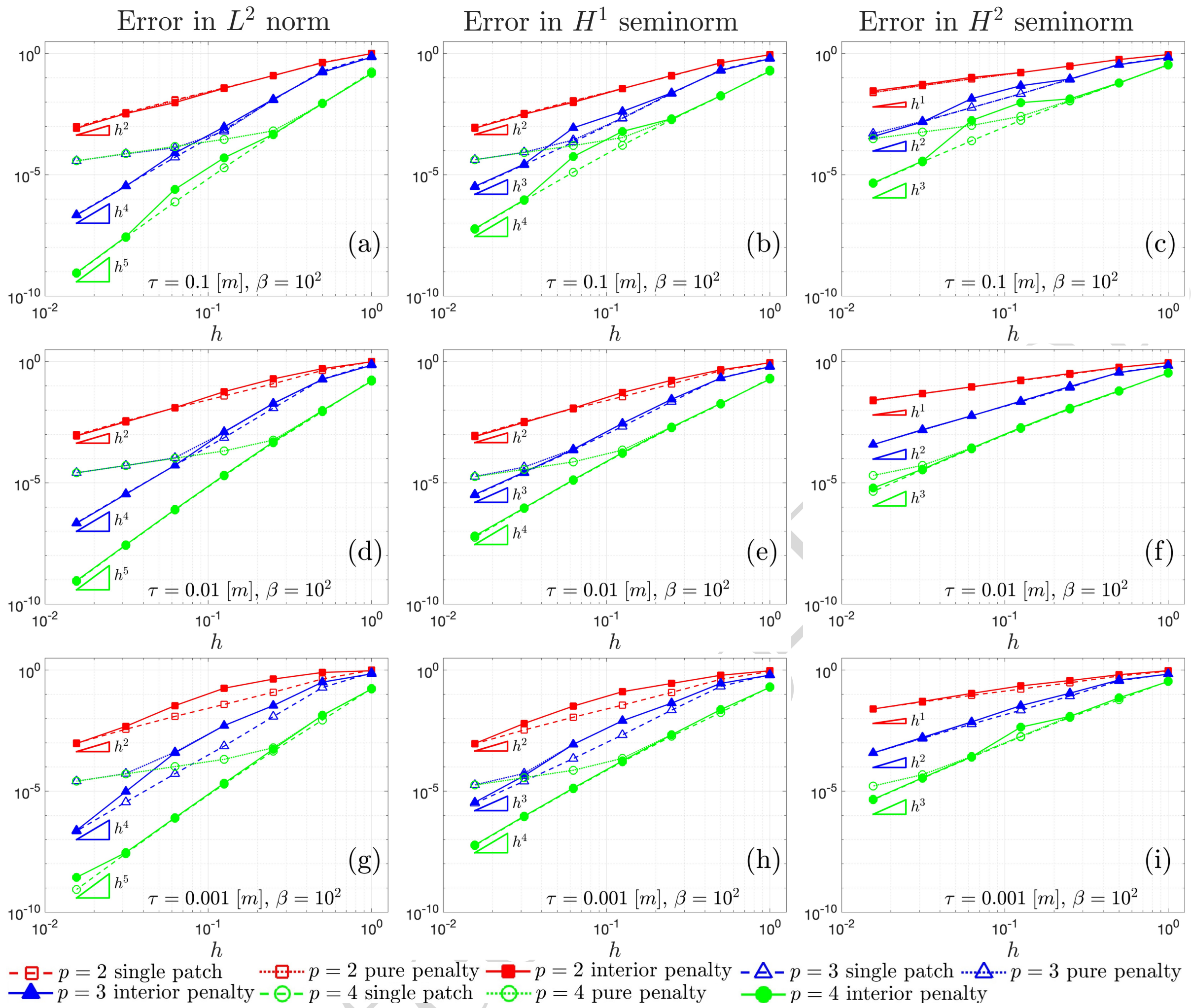

Figure 3: $L^2$ convergence (a), (d) and (g), $H^1$ convergence (b), (e), and (h), and $H^2$ convergence (c), (f) and (i), associated with the Kirchhoff plate shown in Sec.5.1, for $\beta = 10^2$ and for different values of the thickness $\tau$.

Consistently with the findings in [47], the pure penalty method fails to ensure optimal convergence due to the choice of the penalty parameters. Achieving optimal convergence with the pure penalty method would require super-penalization, where the penalty parameters scales with powers of the mesh size that depends on the order of the polynomials [49]. However, such high penalty values typically lead to ill-conditioning of the linear system, especially for high-order polynomials. This underscores the advantage of the interior penalty method, which achieves accurate convergence without facing severe ill-conditioning issues and therefore ensuring accurate results without compromising numerical stability, although some locking phenomena can still be observed in the curves for $p = 2$ and $p = 3$ when $\tau = 0.001$ [m]. The triangles in the graphs show the optimal convergence rates expected for the Kirchhoff-Love theory [60]. It is worth noting that, consistently with the expected theoretical prediction [60, 81], in our results the optimal convergence rate in $L^2$ norm for $p = 2$ is equal to $p$ and not $p + 1$.

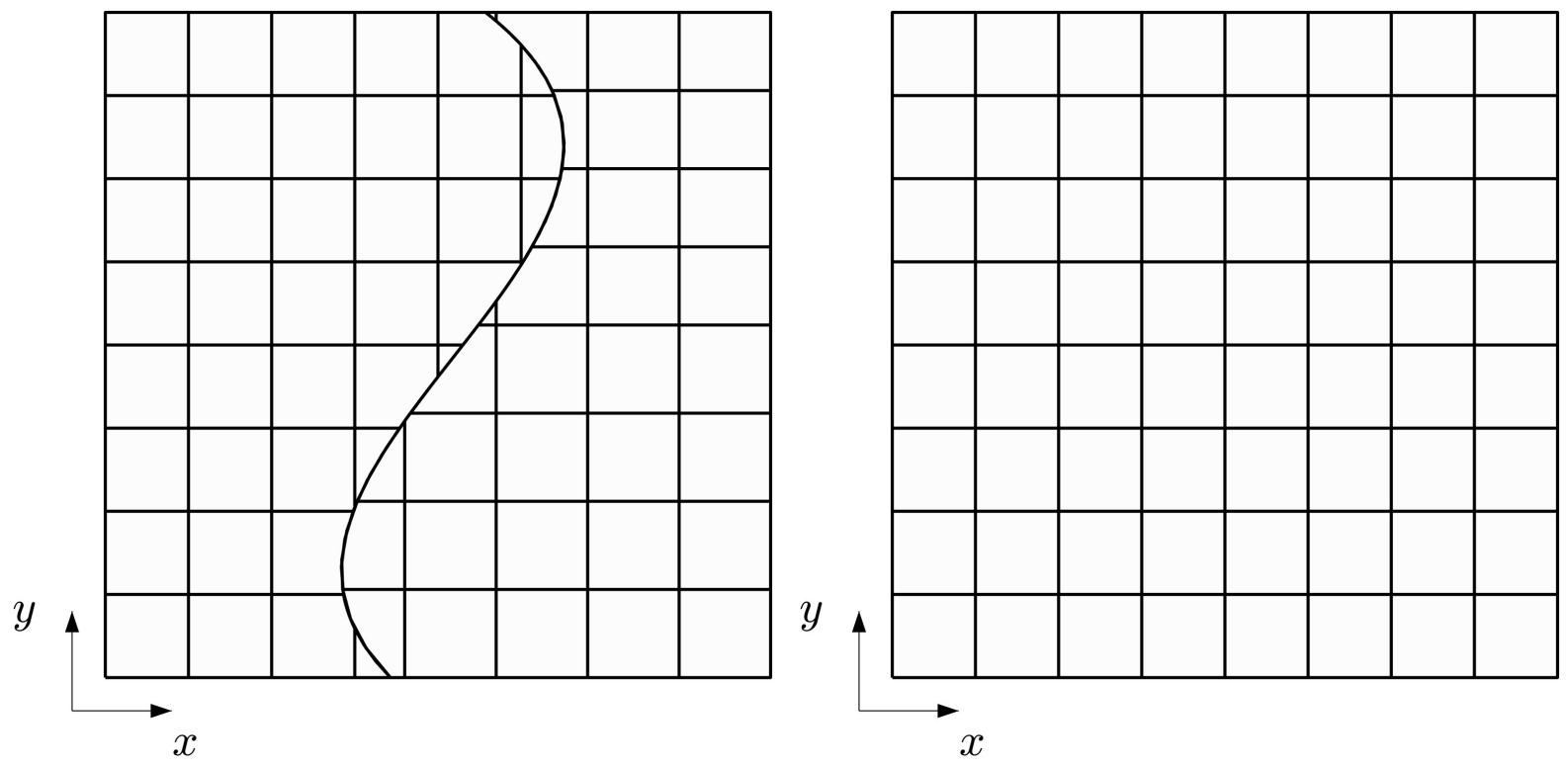

(a) (b)

Figure 4: Discretization of the Kirchhoff laminate in a non-conformal multi-patch (a) and a single-patch (b) setting.

The geometry depicted in Fig.2a is also utilized for conducting a laminate test, with its multi-patch and single-patch configurations presented in Fig.4. In the multi-patch discretization, a knot is inserted in the position $(0.5, 0.5)$ for the left patch and $(0.45, 0.53)$ for the right one. The laminate is constructed using orthotropic laminae with the following properties: longitudinal Young's modulus $E_l^{\langle \ell \rangle} = 25$ [GPa], transversal Young's modulus $E_t^{\langle \ell \rangle} = 1$ [GPa], Poisson's ratio $\nu_{lt}^{\langle \ell \rangle} = 0.25$, shear modulus $G_{lt}^{\langle \ell \rangle} = 0.4$ [GPa] and thickness $\tau^{\langle \ell \rangle} = 0.0025$ [m]. The lamination sequence employed is $[0, 90, 90, 0]$. The boundary conditions are identical to those of the isotropic case, while the force applied to the shell surface is modified in order to manufacture the distribution of the displacement in Eq.(30) with the different material properties. The arbitrary coefficient in the definition of the penalty parameter is selected as $\beta = 10$. The $L^2$, $H^1$, and $H^2$ convergence curves for this test are depicted in Fig.5. Similar to the isotropic case, the observations regarding the convergence properties of the methods investigated hold true also in this scenario.

As a final remark, Fig.6 depicts the laminate's bent structure. The mesh lines, the contour of the displacement vector's magnitude, the components of the generalized moment $M^{11}$, and $M^{12}$, are presented superimposed on the deformed surface of the plate. The images refer to $p = 3$ and $\beta = 10$. It is important to mention that some of the lines in the image are merely used to visualise trimmed elements and do not delimit any actual element edge.

## 5.2 Hyperbolic paraboloid Kirchhoff-Love shell

The second set of numerical experiments focuses on a curved isotropic shell with a mid-surface represented by a hyperbolic paraboloid. This example is derived from the shell obstacle course introduced in [60], which offers a collection of tests having analytical solutions. These tests serve as an appropriate benchmark to assess the performance of the proposed coupling method by evaluating the convergence curve based on the $L^2$ norm, $H^1$, and $H^2$ seminorms of the solution error.

In this test, three values of the shell thickness are taken into account $\tau = 0.1$ [m], $\tau = 0.01$ [m], and $\tau = 0.001$ [m], while the complete definition of the mid-surface can be found in [60]. Nonetheless, for a better understanding, the geometry of the shell and the coordinates of the vertices are illustrated in Fig.7a. The material is isotropic with a Young's modulus of $E = 70$ [GPa] and a Poisson ratio of $\nu = 0.3$. The shell is subjected to non-homogeneous Dirichlet boundary conditions that for the displacement are enforced strongly assigning values directly to the degrees of freedom after the $L^2$ projection of the displacement field in the spline space of the corresponding edges. Conversely, the bending rotation corresponding to the reference displacement field is weakly imposed at the boundary using the Nitsche's method described in Sec.4.2, using in Eqs.(21) and (23) only the terms corresponding to the bending rotation. The boundary conditions and the domain force are chosen to manufacture the reference displacement field:

$$\boldsymbol{u}^{ex} = U_0 \xi_2 \sin\left(\frac{\pi}{2}\xi_2\right)\boldsymbol{e}_1 + U_0 \xi_2 \sin\left(\frac{\pi}{2}\xi_2\right)\boldsymbol{e}_2, \tag{31}$$

where $U_0 = 1$ [m]. The function for the force is provided in the Mathematica notebook associated with [60]. Similarly to the previous test, the shell is modeled using a two-patch configuration (see Fig.7b), as well as a single patch configuration (see Fig.7c) which serves as a reference for the convergence curves. To make the two patches in the multi-patch configuration non-conforming at the interface, a knot is inserted in each patch at the curvilinear coordinates $(0.5, 0.5)$ and $(0.54, 0.43)$, respectively. In such way, even after subsequent dyadic refinements of the discretization, the non-conforming nature of the interface is maintained.

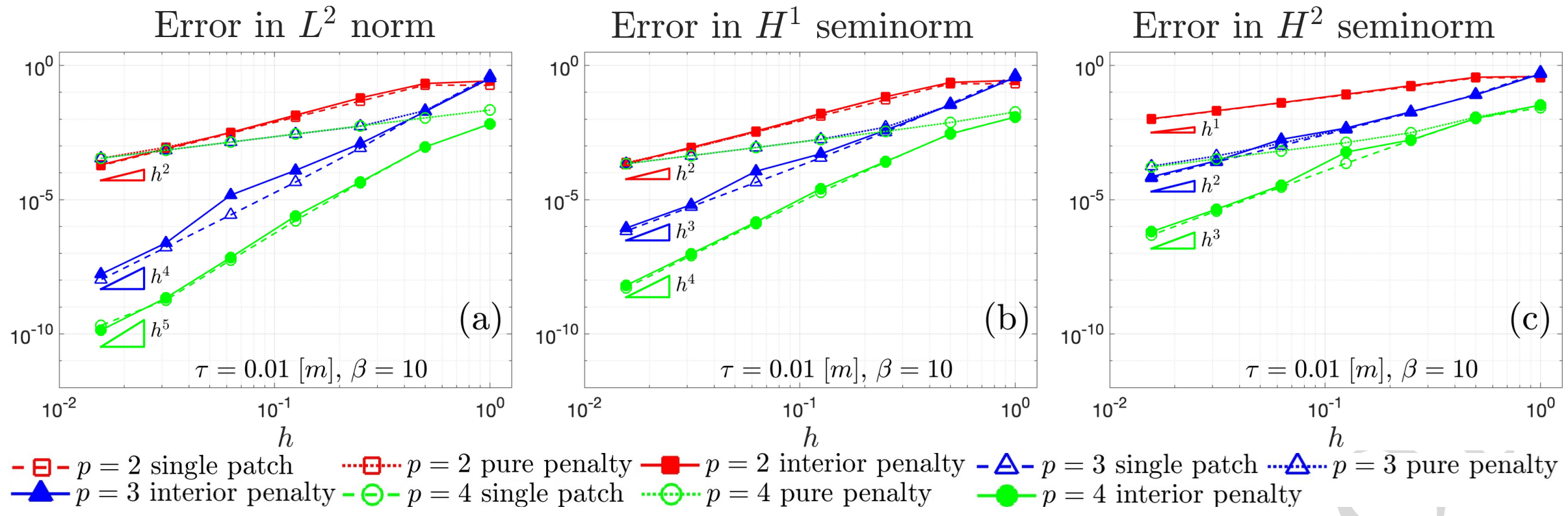

Figure 5: $L^2$ norm convergence (a), $H^1$ convergence (b), and $H^2$ convergence (c), correspondent to the Kirchhoff laminate having thickness $\tau = 0.01$ [m].

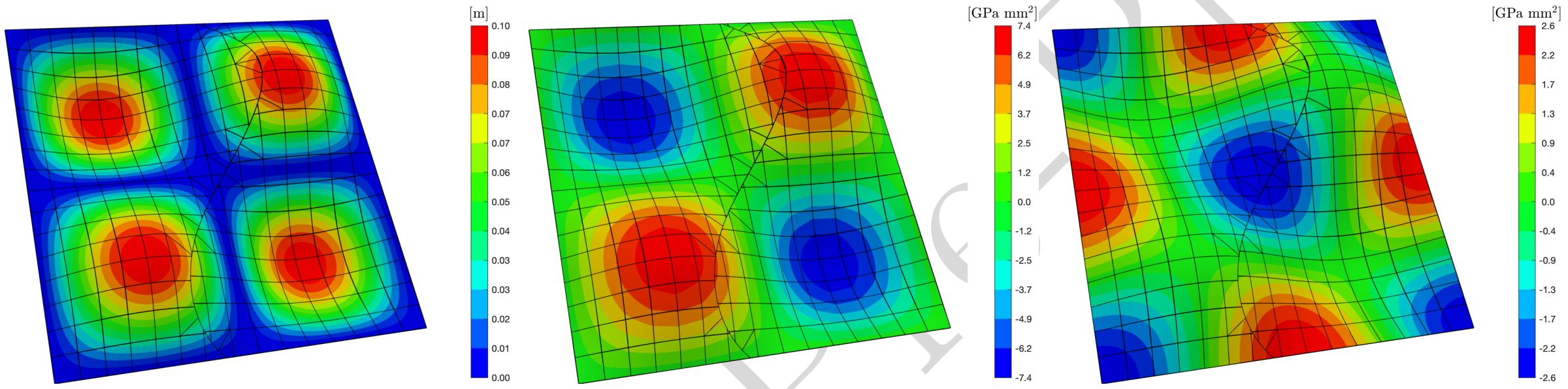

(a) $|\boldsymbol{u}|$ (b) $M^{11}$ (c) $M^{12}$

Figure 6: Deformed configuration for the laminate in Fig.4a with superimposed mesh grid and contour of the magnitude of the displacement (a), generalized moments $M^{11}$ (b), and $M^{12}$ (c).

Figs. 8 and 9 show the $L^2$, $H^1$, and $H^2$ convergence curves for the proposed test, considering different values of the shell thickness, the polynomial order, and the approximation approach. Two values of the arbitrary parameter $\beta$ that appears in the penalty terms in Eq.(29) are taken into account, specifically, $\beta = 10$ in Fig.8 and $\beta = 10^3$ in Fig.9. Three different approaches are compared: a single-patch, a two-patch configuration coupled through a pure penalty method, and a two-patch configuration coupled through the interior penalty method.

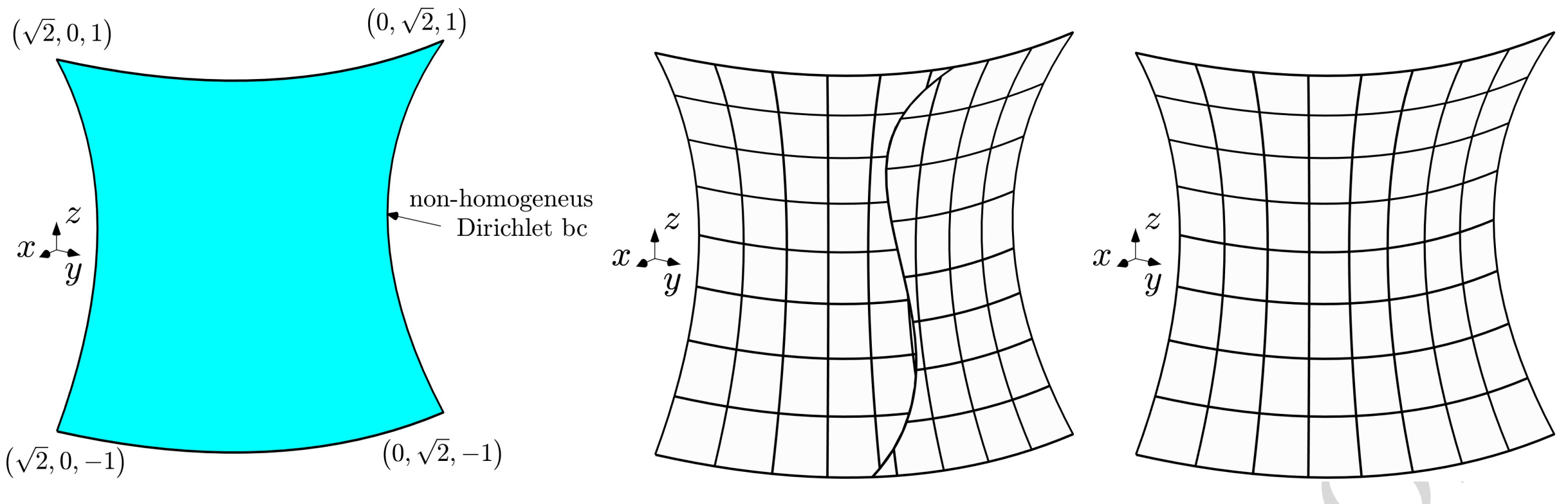

Figure 7: Geometry of the hyperbolic paraboloid described in Sec.5.2 (a). Discretization of the shell employing two non-conformal IGA patches (b) and a single IGA patch (c), shown for a certain refinement level.

The first feature that can be noticed in both Figs. 8 and 9 is that, with regard to the $L^2$ norm, the interior penalty method tend to closely follow the convergence curves of the reference single patch discretization, in contrast to the pure penalty method. In this comparison, the penalty value is scaled in the same way for both the interior penalty and pure penalty methods. Once again, it is worth remarking that achieving optimal convergence rates with the pure penalty approach would necessitate a scaling of the penalty parameter with higher exponents of the mesh size $h$, thereby directly deteriorating the condition number of the linear system.

Focusing on $L^2$ convergence for $p = 4$ and $\tau = 0.1$ [m], the curve for the single patch undergoes a slight change with an increase in the parameter $\beta$. This change occurs because, while the boundary conditions for displacements are strongly imposed, the boundary conditions for rotation are still enforced through Nitsche's method. Consequently, the condition number for the last refinement level is influenced by the higher penalty, causing the curve to deviate from the optimal rate. This phenomenon is also evident in the interior penalty method, where weak coupling among trimmed elements also adversely affects the condition number. In the curves corresponding to $L^2$ norm for $p = 4$, convergence is lost in the last refinement level for every combination of $\tau$ and $\beta$.

A similar behavior is observed in $H^1$ convergence, although the spurious increase of the error in the last refinement level is more pronounced. This is particularly noticeable in the graphs associated with $\tau = 0.1$ [m] and $\beta = 10$, where also the curves for $p = 2$ and $p = 3$ are affected. In the latter case, it is apparent that further enhancing the stability of the interior penalty method by increasing the arbitrary parameter to $\beta = 10^3$ effectively reduces the error. This observation suggests that, in certain scenarios, stability issues may arise independently of the condition number of the linear system.

An in-depth examination is essential for the graphs in $H^2$ seminorm. Notably, in these graphs, spurious effects become evident even in earlier refinement levels, and a clear asymptotic regime, demonstrating optimal convergence, is achieved only in certain instances. The spurious increase in error in the last refinement level diminish only slightly with an increase in the penalty parameter. However, it is apparent that, for $\beta = 10^3$, the earlier refinement levels closely concord to the reference single patch discretization, suggesting that a more robust stabilization technique would indeed be beneficial. This test confirms that, in the $H^2$ seminorm for Kirchhoff-Love trimmed shell patches coupled weakly, criticalities tend more easily to emerge.

Another noteworthy observation is that higher values of $\beta$ result in increased errors in the preasymp-totic regime for the interior penalty method. This effect is expected since the non-conforming nature of the coupling interface leads to a locking effect due to the stabilization terms in the integrals in Eq.(27b). In fact, the discrete approximation spaces of the displacement fields for the two patches are unable to perfectly match at the interface for non-trivial distribution, causing spurious locking phenomena. Further increasing the penalty exacerbates this issue.

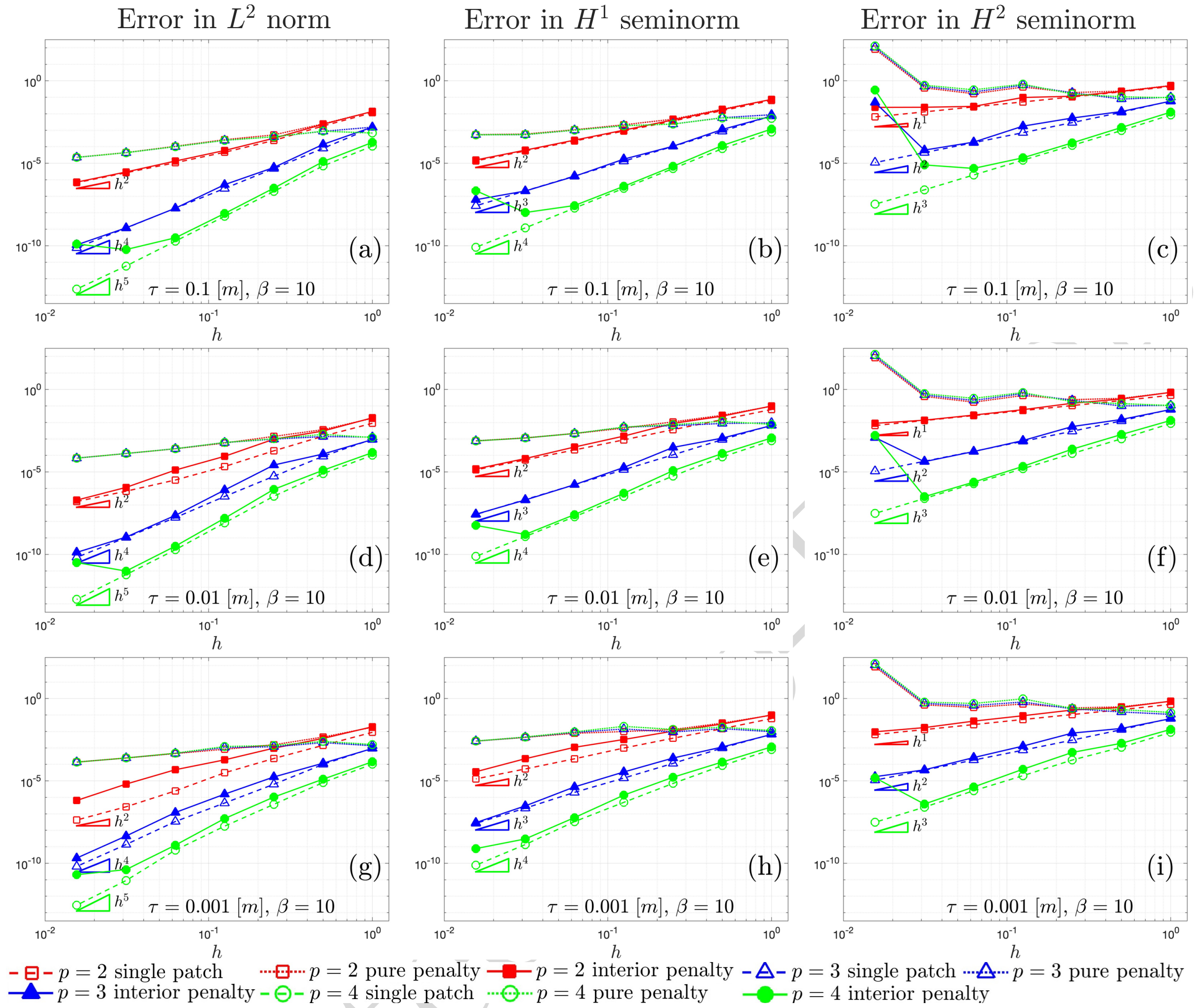

Figure 8: $L^2$ convergence (a), (d) and (g), $H^1$ convergence (b), (e), and (h), and $H^2$ convergence (c), (f) and (i), associated with the shell shown in Fig.7a, for different values of the thickness $\tau$ and arbitrary parameter $\beta = 10$.

To better assess the stability of the method, the discretization shown in Fig.11 is adopted. The configuration involves two patches with inserted knots at $(0.5-\delta,\ 0.5)$ and $(0.5+\delta, 0.43)$, respectively. Both patches are trimmed along the line characterized by $\xi_1 = 0.5$, and then weakly merged together. This arrangement ensures that the coupling interface leans on two columns of arbitrarily small trimmed elements. Each patch consists of a $16 \times 9$ grid with $p = 3$, with one column of 16 elements being critically cut. The value of $\delta$ is varied, and the trend on the error as $\delta$ approaches zero is depicted in Fig.12. The increase in error is evident in all the norms considered, and more pronounced in $H^2$ seminorm, as expected. The graph in Fig.12d shows the value of the penalty parameter and the condition number after the application of the Jacobi preconditioner, showing that the proposed stabilization technique still affect the condition number. This test demonstrates that when dealing with trimmed elements, in certain critical conditions, the method can indeed lead to instabilities. Addressing this issue is a priority in future developments of the present work to increase robustness.

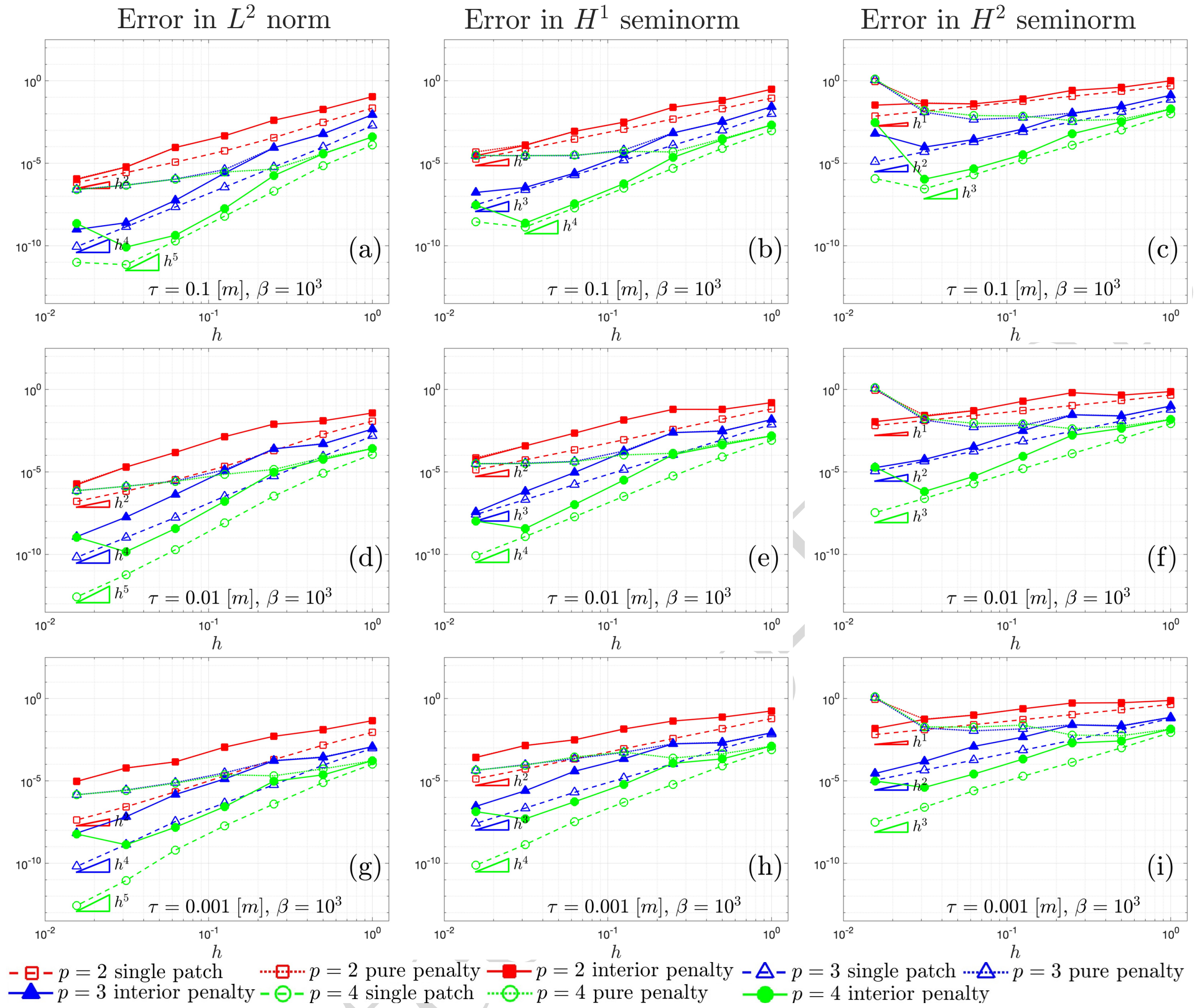

Figure 9: $L^2$ convergence (a), (d) and (g), $H^1$ convergence (b), (e), and (h), and $H^2$ convergence (c), (f) and (i), associated with the shell shown in Fig.7a, for different values of the thickness $\tau$ and arbitrary parameter $\beta = 10^3$.

It is also worth noting that the proposed method is adopted for a wide range of thickness values. However, it is important to highlight that the Kirchhoff-Love shell equation may not be suitable for thickness ratios that are too high. In such cases, higher-order theories should be employed to accurately capture the behavior of the shells.

Lastly, for the sake of completeness, in Fig.10 it is shown the contour of the magnitude of the displacement vector, the contour of the generalized force component $N^{11}$, and the contour of the generalized moment component $M^{11}$, superimposed to the undeformed shell mid-surface, together with the mesh edges. The images refer to $p = 4$ and $\beta = 10^3$.

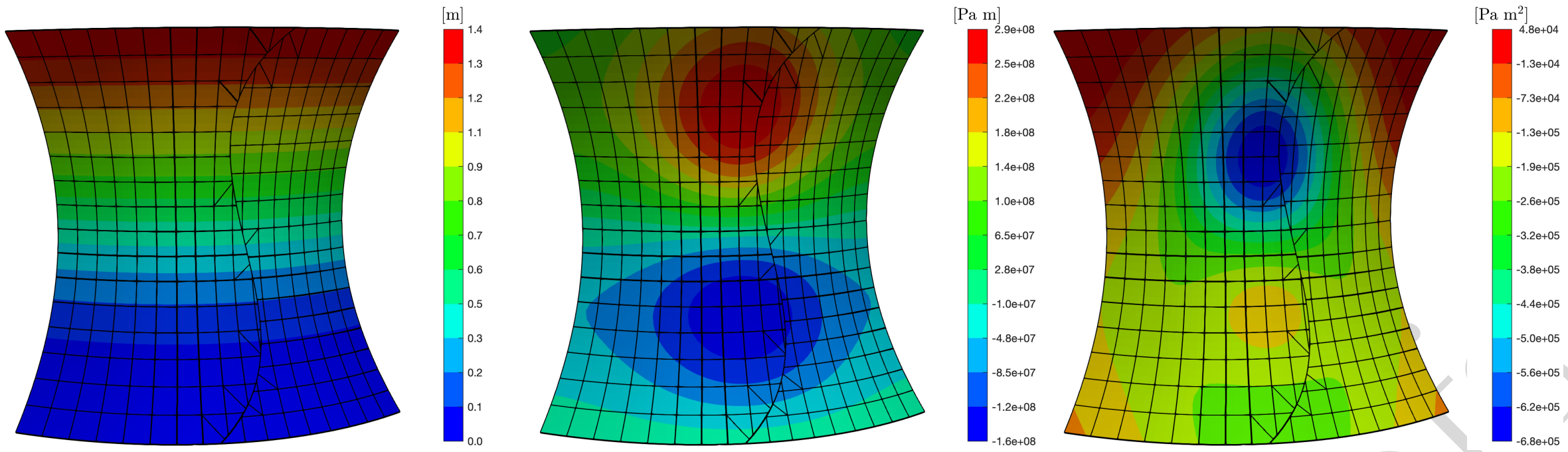

(a) $|\boldsymbol{u}|$ (b) $N^{11}$ (c) $M^{11}$

Figure 10: Undeformed configuration for the shell in Fig.7a with superimposed mesh grid and contour of the magnitude of the displacement (a), generalized force $N^{11}$ (b), and generalized moment $M^{11}$.

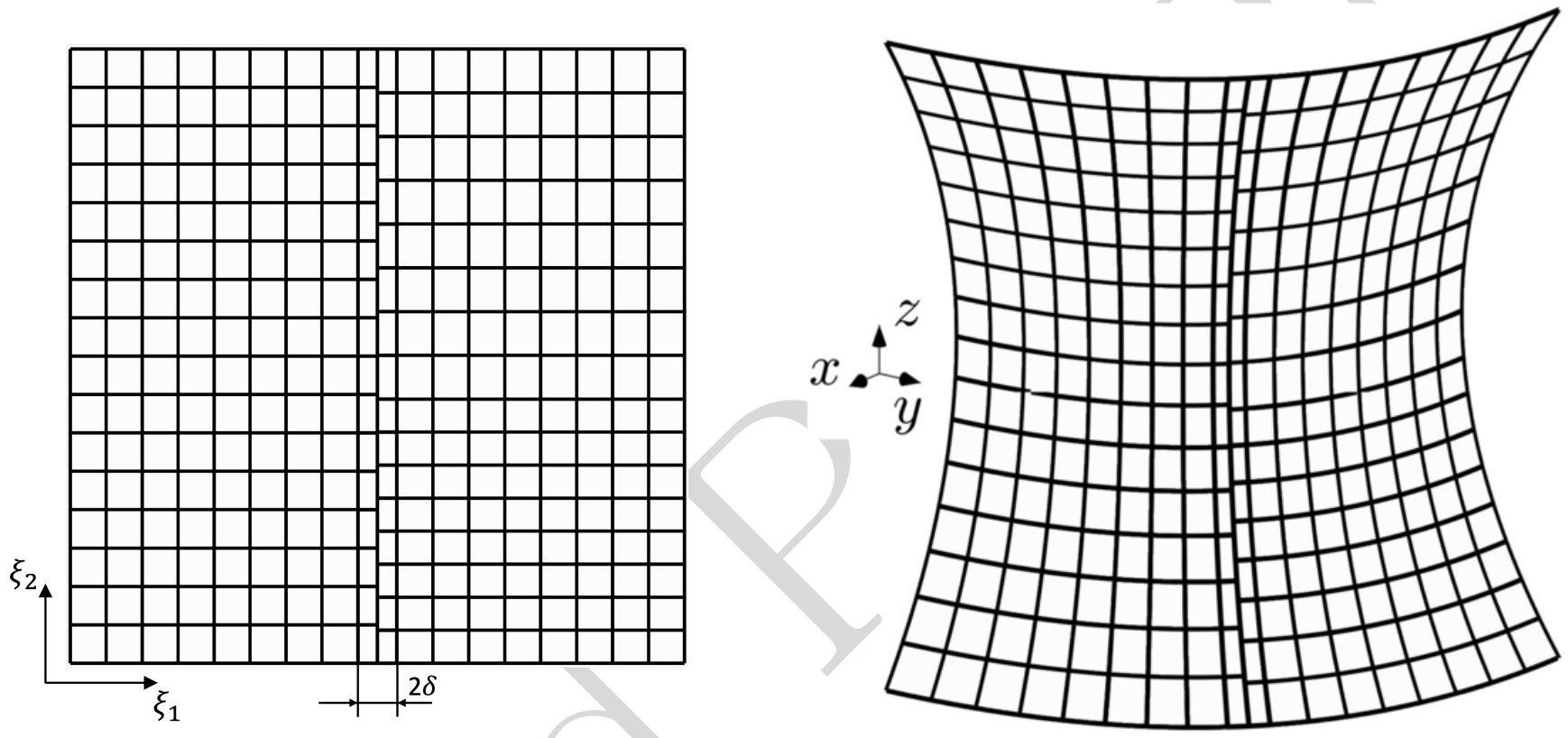

(a) (b)

Figure 11: Discretization of the hyperbolic paraboloid shell in the parametric domain (a) and in the Euclidean space (b), consisting of two patches meeting at a non-conforming straight interface that leans on two columns of trimmed and poorly shaped elements.

## 5.3 Coupling of intersecting cylindrical shells

The final test aims to demonstrate the efficiency and robustness of the proposed method, as well as the effectiveness of the trimming and coupling algorithms. This test involves a geometry with a complexity level comparable to that encountered in real-world industrial applications, consisting of multiple patches intersecting at variable angles. The structure under investigation consists of five intersecting cylinders, as depicted in Fig.13a. The main cylinder has a length of $L = 8$ [m] and a radius of $R = 1$ [m], while the remaining cylinders have an untrimmed length of $L = 4$ [m] and a radius of $R = 0.8$ [m]. The complete geometrical description of the structure is not provided here for the sake of conciseness but can be found in the STEP file associated with this publication.

Figs.13b to 13f illustrate the trimmed parametric domains of the corresponding cylindrical patches in Fig.13a. It is noteworthy that the intersecting curves, both in the physical and parametric domains, consist of multiple segments connected in general with $C^0$ continuity. For such geometric configurations, the robustness of the algorithm for identifying quadrature points with high-order precision becomes crucial.

The cylinders in the structure are simply-supported at the external edges, and the structure is subjected to a uniformly applied domain traction given by $\tilde{\boldsymbol{F}} = [10^5, 10^5, 10^5]^T$ [Pa]. The material used is a laminate with layers made of the same orthotropic material as described in Sec.5.1, with a thickness of $\tau^{\langle \ell \rangle} = 0.0025$

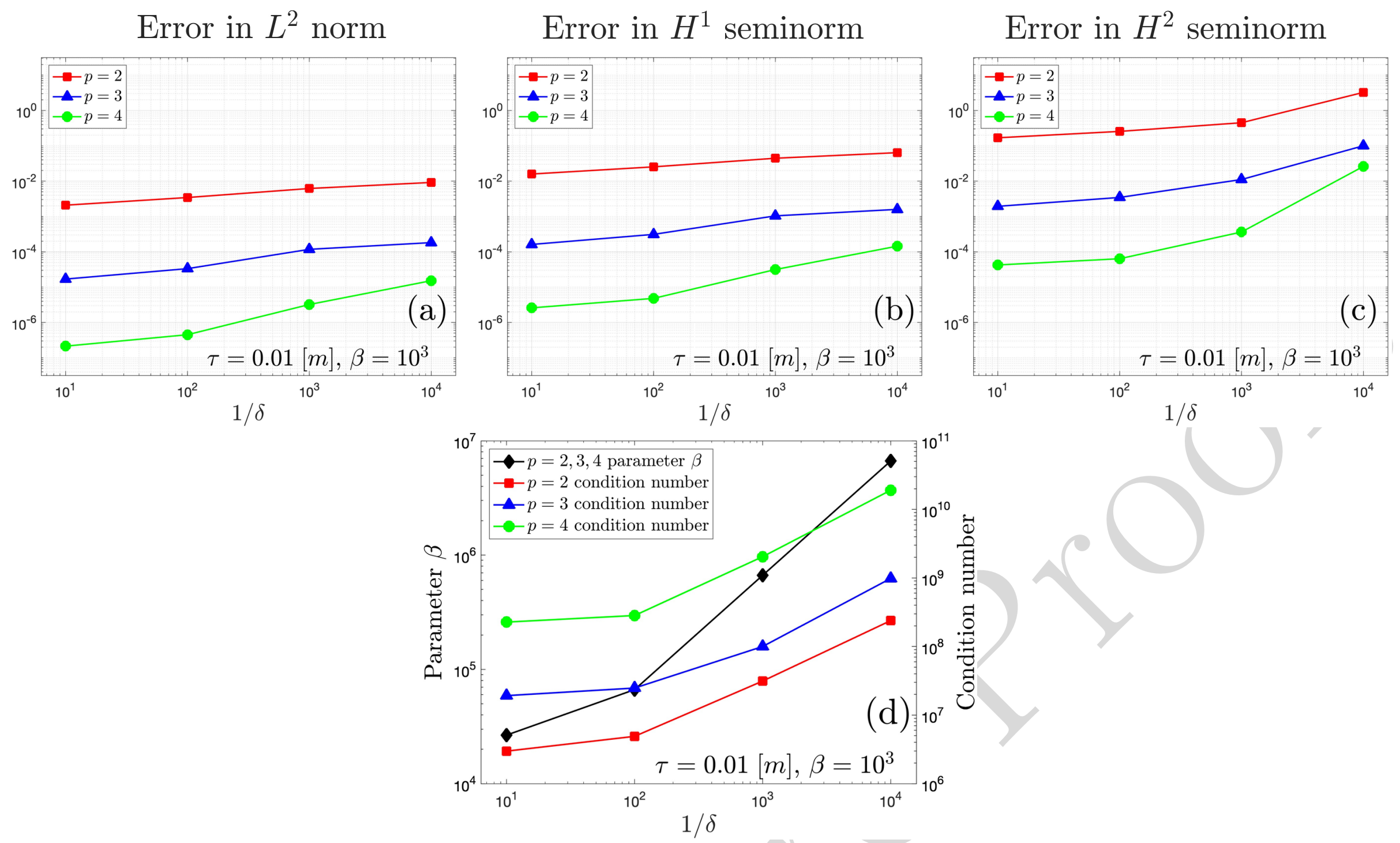

Figure 12: Error on $L^2$ norm (a), $H^1$ seminorm (b), and $H^2$ seminorm (c), value of the parameter $\beta$ (left ordinate in (d)), and condition number (right ordinate in (d)), as functions of decreasing value of the parameter $\delta$ for the discretization shown in Fig.11.

[m], and lamination sequence $[90, 0, 0, 90]$;

In each patch, the continuity along the circumferential directions between the edges corresponding to the first and last values of the knot vector is enforced here using periodic boundary conditions obtained adopting periodic spline spaces [74]. It is worth noting that along these edges, a weak imposition of the coupling condition could also be applied by considering an interface with both edges coming from the same patch. However, using periodic boundary conditions allows for a reduction in the overall number of degrees of freedom in the analysis. The cylindrical patches meet at interfaces with non-zero angles, that are easily managed by the proposed formulation, as it is not limited to $G^1$ surfaces. The polynomial degree used for each patch in each direction is $p = 6$, that can be easily adopted thanks to the straightforward construction of B-spline basis functions. It is worth mentioning that to properly integrate in the trimmed elements and their respective boundaries nine Gaussian points were adopted in each direction.

Fig.14 presents the contour plot of the magnitude of the displacement from two different views. The results obtained using the formulation described in this paper (a) and (b) are compared with those obtained using triangular elements (STRI3) in the Abaqus® software [92], (c) and (d). The close agreement between the two approaches highlights the competitiveness of the proposed method with the available finite element software, affirming its accuracy and reliability even for complex geometries consisting of patches intersecting at a variable angle. For completeness, in Fig.15 the contours of the components of the generalized force $N^{22}$ and $N^{12}$, as well as the components of the generalized moment $M^{11}$ and $M^{22}$, as defined in Eq.(12) are also depicted.

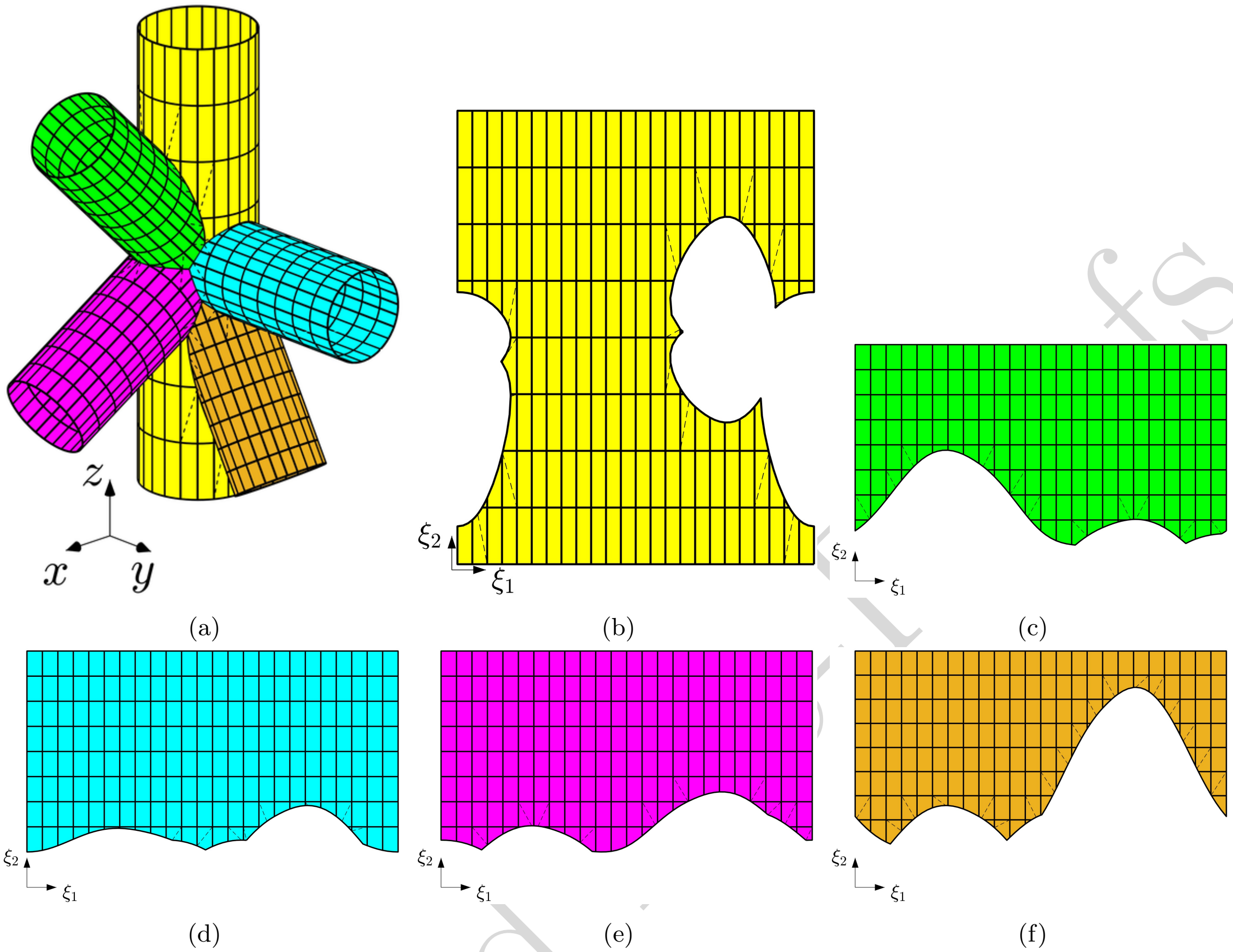

Figure 13: Geometry of the test described in Sec.5.3 with superimposed mesh (a). Parametric domain of each of the B-rep patches (b – f).

# 6 Conclusions

In this work, we have conducted the linear elastic static analysis of isotropic and laminated plates and shells using the Kirchhoff-Love theory. The associated fourth-order problem imposes a $C^1$ continuity requirement on the solution space that through the Isogeometric Analysis (IGA) was seamlessly addressed by simply using high-degree splines.

The proposed method is capable of handling structures formed by multiple IGA patches meeting at interfaces that may not necessarily be conforming and might even intersect at an angle. In fact, amongst the various coupling strategies, we employed the symmetric interior penalty method. This method produces a weak imposition of the coupling conditions, that therefore does not need to be embedded in the solution space, easily lending itself to non-conforming and trimmed intersections amongst patches. Unlike other Nitsche-type methods, the interior penalty method yields a symmetric solving linear system. The proposed choice for the penalty parameters in the stabilization terms of the weak variational statements demonstrated to be effective in the numerical results, where multi-patch discretization with non conforming interfaces were tested. The error norms showed optimal convergence in the asymptotic regime, in agreement with expected theoretical rates, and proved to be competitive respect to single patch discretizations. For comparison, results related to a pure penalty formulation were also shown but exhibited a behaviour not comparable with the interior penalty's one.

In the proposed formulation, Dirichlet boundary conditions are applied both strongly and weakly,

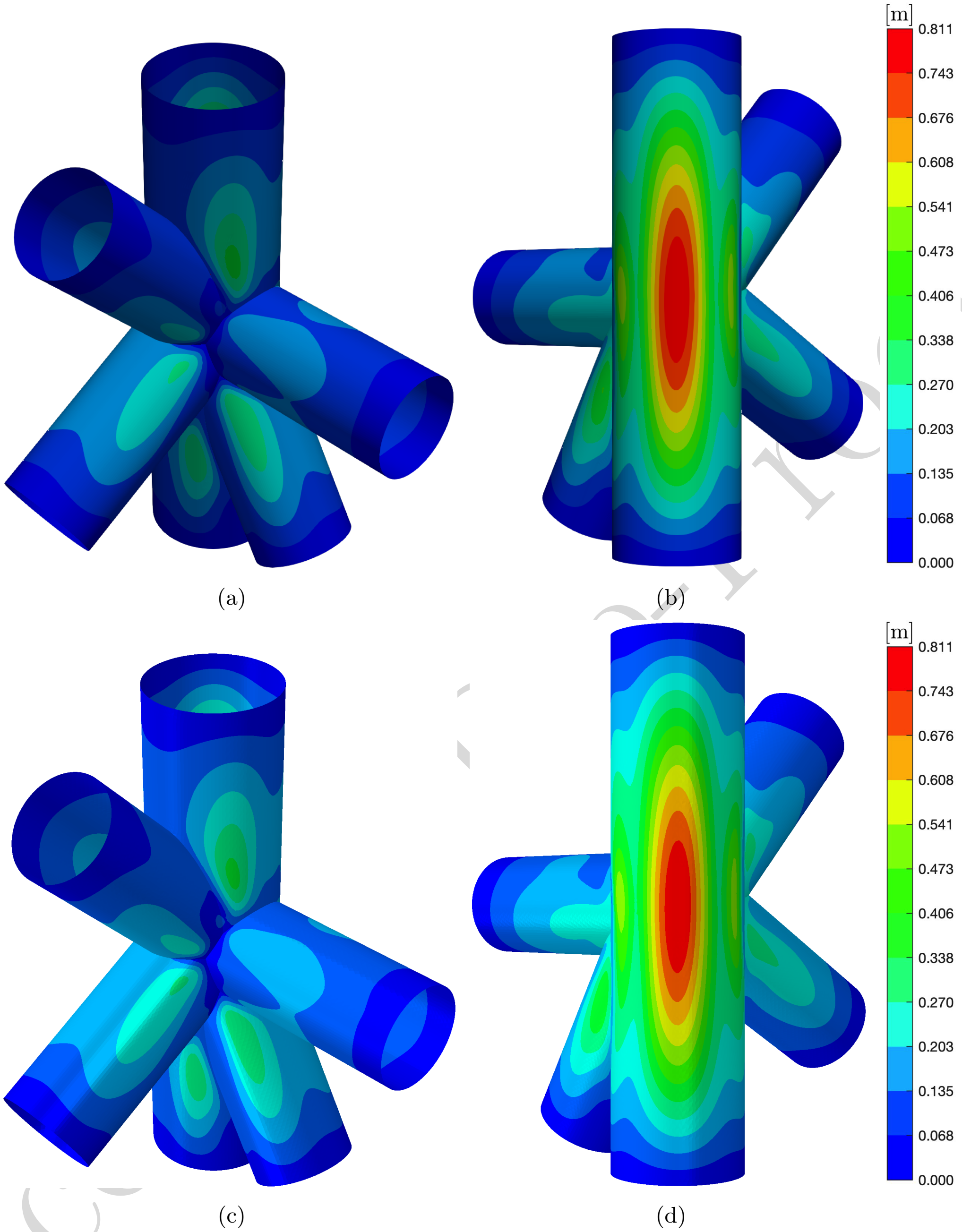

Figure 14: Two different views of the contour of the magnitude of the displacement for the structure described in Sec.5.3 obtained with the method presented here (a) and (b), and with the elements STRI3 in Abaqus® (c) and (d).

depending on whether the boundary is conforming or trimmed, and whether the Dirichlet condition refers to displacement or rotation. The expression for the fluxes used to apply coupling conditions and essential boundary conditions was adopted from [60], where their correctness was rigorously demonstrated, in contrast with the formulation typically used in the existing literature. This work successfully replicated the results for the hyperbolic paraboloid benchmark from the new shell obstacle course, for the first time in the context of non-conforming multi-patch discretizations. The method proved optimal convergence, although for high refinement level and for high polynomials some ill-conditioning of the linear system appeared.

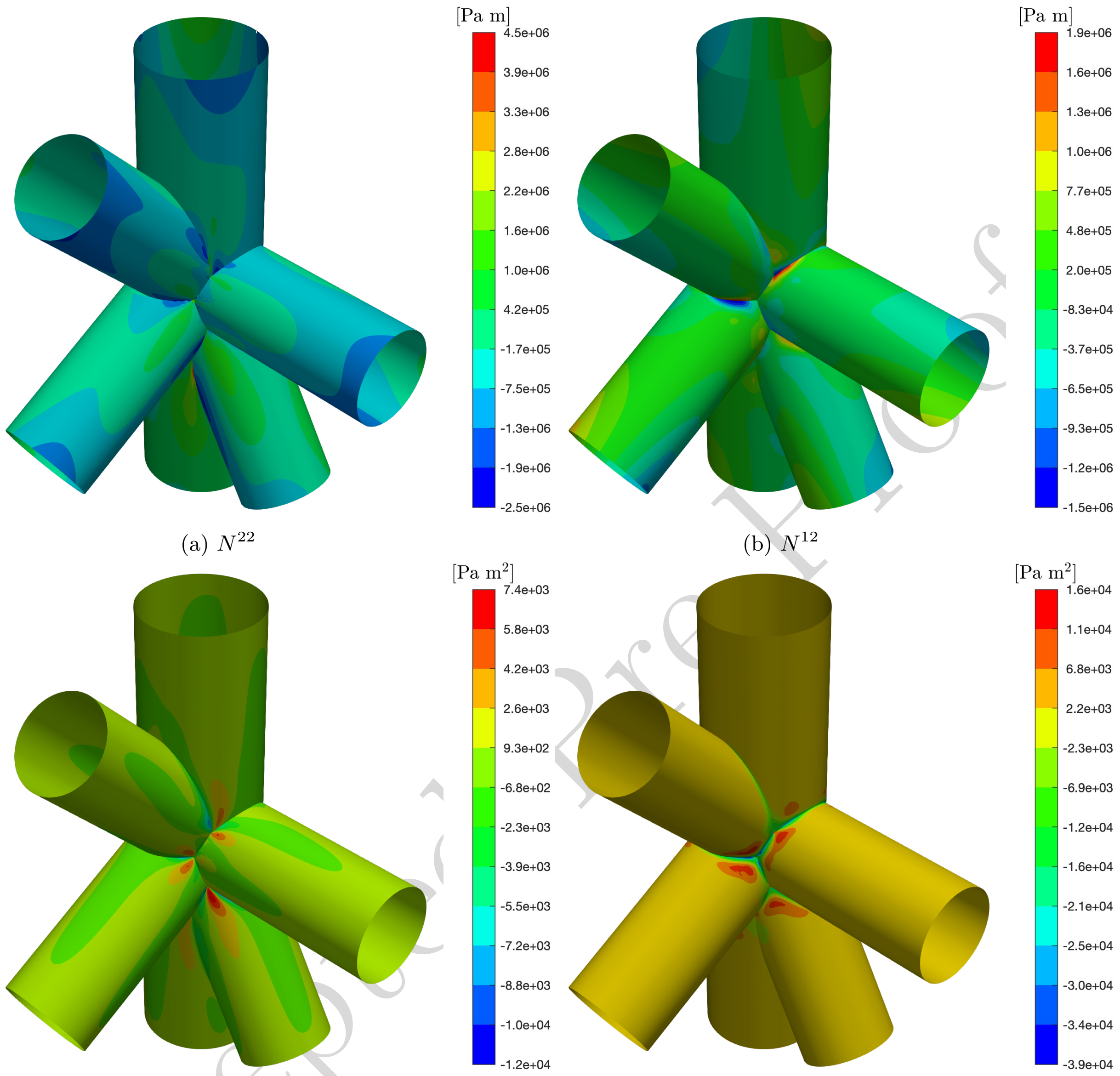

(a) $N^{22}$ (b) $N^{12}$

(c) $M^{11}$ (d) $M^{22}$

Figure 15: Contours of some representative components of the generalized force and moment for the test regarding the intersecting cylindrical shells.

The presented formulation allows for laminated shells, assuming uniform thickness and lamination angle for the layers. However, the extension to variable thickness or variable angle tow composite is straight-forward. The numerical experiments on laminates include a laminated Kirchhoff plate and a structure comprising five laminated intersecting cylinders. For the first one, an analytical solution exists, and the error convergence curves were shown, confirming the properties of the method already demonstrated for an isotropic material. In the last test, a comparison of the solution contour with that obtained with Abaqus® was shown. This last test showcases the potential of the proposed method to handle laminated shell structures of industrial application level of complexity. The software adopted in this work is in fact capable of importing directly STEP files, making it a powerful tool for mechanical analysis.

To what concerns future development, possible directions for further research include: extending the method to address other mechanical phenomena, such as linear buckling, free vibration, and transient response, as well as nonlinear elasticity; investigating the limits of the interior penalty method and developing efficient approaches to address potential loss of stability as well as the ill-conditioning issues, especially for trimmed patches; extending the method to encompass Reissner-Mindlin and higher-order shell theories.

# Acknowledgements

This research was supported by the European Union Horizon 2020 research and innovation program, under grant agreement No 862025 (ADAM2), and by the Swiss National Science Foundation through the project No 200021_214987 (FLASh). G. Guarino further acknowledges the support of the Swiss Government Excellence Scholarships program by the Federal Commission for Scholarships for Foreign Students through the scholarship No 2022.0014.

This version of the article has been accepted for publication, after peer review but is not the Version of Record and does not reflect post-acceptance improvements, or any corrections. The Version of Record is available online at: https://doi.org/10.1007/s00366-024-01965-5.

# Declaration of competing interest

The authors declare that they have no known competing financial interests or personal relationships that could have appeared to influence the work reported in this paper.

# A Linear membrane and bending strain and rotation

In this Appendix the problem formulated in Sec.3 is completed by providing the definition of $\boldsymbol{\varepsilon}$, $\boldsymbol{\kappa}$, and $\theta_n$, which are obtained by linearizing their non-linear expressions, at the condition of uniform zero displacement. For a more comprehensive understanding of the underlying continuum mechanics concepts, the interest readers are referred to [59]. The expressions for the components of linear membrane and bending strain are:

$$\varepsilon_{\alpha\beta} = \frac{1}{2}(\boldsymbol{a}_\alpha \cdot \boldsymbol{u}_{,\beta} + \boldsymbol{a}_\beta \cdot \boldsymbol{u}_{,\alpha}) \,, \tag{32a}$$

$$\kappa_{\alpha\beta} = -\boldsymbol{a}_3 \cdot \boldsymbol{u}_{,\alpha\beta} - \left(\frac{\boldsymbol{a}_{\alpha,\beta} \times \boldsymbol{a}_1}{\lambda} + \boldsymbol{a}_{\alpha,\beta}^T \boldsymbol{M}\boldsymbol{M}_1\right) \cdot \boldsymbol{u}_{,2} + \left(\frac{\boldsymbol{a}_{\alpha,\beta} \times \boldsymbol{a}_2}{\lambda} + \boldsymbol{a}_{\alpha,\beta}^T \boldsymbol{M}\boldsymbol{M}_2\right) \cdot \boldsymbol{u}_{,1} \,. \tag{32b}$$

In Eq.(32) the matrices $\boldsymbol{M}$ and $\boldsymbol{M}_\alpha$ have been introduced as

$$\boldsymbol{M} = \frac{\boldsymbol{a}_3 \otimes \boldsymbol{a}_3}{\lambda} \,, \tag{33a}$$

$$\boldsymbol{M}_\alpha = [\boldsymbol{e}_1 \times \boldsymbol{a}_\alpha \quad \boldsymbol{e}_2 \times \boldsymbol{a}_\alpha \quad \boldsymbol{e}_3 \times \boldsymbol{a}_\alpha] \,, \tag{33b}$$

where $\lambda = |\boldsymbol{a}_1 \times \boldsymbol{a}_2|$. Similarly, the linear expression for the normal rotation is obtained from the non-linear one in [45] obtaining

$$\theta_n = \left(\frac{(\boldsymbol{t} \times \boldsymbol{a}_3) \times \boldsymbol{a}_1}{\lambda} + (\boldsymbol{t} \times \boldsymbol{a}_3)^T \boldsymbol{M}\boldsymbol{M}_1\right) \cdot \boldsymbol{u}_{,2} - \left(\frac{(\boldsymbol{t} \times \boldsymbol{a}_3) \times \boldsymbol{a}_2}{\lambda} + (\boldsymbol{t} \times \boldsymbol{a}_3)^T \boldsymbol{M}\boldsymbol{M}_2\right) \cdot \boldsymbol{u}_{,1} \,. \tag{34}$$

It is remarked that the expressions shown in Eq.(32) for the strain are equivalent to those in [27]. However, we find the proposed notation to be more convenient as it allows for the introduction of some terms that are common with the definition of $\theta_n$.

# B Stiffness matrices for a laminate in tensor notation

Due to the assumption of orthotropic homogeneous layers, each of the $\ell$ layers is characterized by three orthotropic directions. We assume that thickness direction $\boldsymbol{a}_3$ is always one of these three. The other two directions, namely the longitudinal and transversal direction respect to the fibers, depend on the lamination angle of the layer $\theta^{\langle\ell\rangle}$. In the orthotropic reference system the stiffness matrix in Voigt notation is

$$\boldsymbol{E}_{ort}^{\langle\ell\rangle} = \begin{bmatrix} 1/E_l^{\langle\ell\rangle} & -\nu_{lt}^{\langle\ell\rangle}/E_l^{\langle\ell\rangle} & 0 \\ -\nu_{tl}^{\langle\ell\rangle}/E_t^{\langle\ell\rangle} & 1/E_t^{\langle\ell\rangle} & 0 \\ 0 & 0 & 1/G_{lt}^{\langle\ell\rangle} \end{bmatrix}^{-1} , \tag{35}$$

where $E_l^{\langle\ell\rangle}$ and $E_t^{\langle\ell\rangle}$ are the Young moduli in the longitudinal and transversal directions, $\nu_{lt}^{\langle\ell\rangle}$ and $\nu_{tl}^{\langle\ell\rangle}$ are the Poisson ratios, and $G_{lt}^{\langle\ell\rangle}$ is the shear modulus, for the $\ell$-th layer. It is worth mentioning that Eq.(35) takes into consideration also the plane stress state hypothesis [93]. To retrieve the case of an isotropic material, a single-layer laminate is considered where it is assumed $E_l^{\langle 1\rangle} = E_t^{\langle 1\rangle} = E$, $\nu_{lt}^{\langle 1\rangle} = \nu_{tl}^{\langle 1\rangle} = \nu$, and $G_{lt}^{\langle 1\rangle} = E/(1+\nu)$, where $E$ and $\nu$ are the only Young modulus and Poisson ratio for the material under investigation.

To define the equivalent stiffness matrices for the entire laminate, the stiffness matrix of each layer is rotated in the common orthonormal reference system $\boldsymbol{m}_1\boldsymbol{m}_2$, which is defined as

$$\boldsymbol{m}_1 = \frac{\boldsymbol{a}_1}{|\boldsymbol{a}_1|} , \tag{36a}$$

$$\boldsymbol{m}_2 = \frac{\boldsymbol{a}^2}{|\boldsymbol{a}^2|} , \tag{36b}$$

through the equation

$$\boldsymbol{E}^{\langle\ell\rangle} = \boldsymbol{T}^{\langle\ell\rangle}\boldsymbol{E}_{ort}^{\langle\ell\rangle}\boldsymbol{T}^{\langle\ell\rangle T} , \tag{37}$$

where the rotation matrix $\boldsymbol{T}$ is the following

$$\boldsymbol{T}^{\langle\ell\rangle} = \begin{bmatrix} \cos^2\theta^{\langle\ell\rangle} & \sin^2\theta^{\langle\ell\rangle} & -2\sin\theta^{\langle\ell\rangle}\cos\theta \\ \sin^2\theta^{\langle\ell\rangle} & \cos^2\theta^{\langle\ell\rangle} & 2\sin\theta^{\langle\ell\rangle}\cos\theta^{\langle\ell\rangle} \\ \sin\theta^{\langle\ell\rangle}\cos\theta^{\langle\ell\rangle} & -\sin\theta^{\langle\ell\rangle}\cos\theta^{\langle\ell\rangle} & \cos^2\theta^{\langle\ell\rangle} - \sin^2\theta^{\langle\ell\rangle} \end{bmatrix} . \tag{38}$$

The laminate generalized stiffness matrices in the reference $\boldsymbol{m}_1\boldsymbol{m}_2$ are obtained after a through-the-thickness integration as:

$$\bar{\boldsymbol{A}} = \sum_{\ell=1}^{N_\ell} \int_{\tau_b^{\langle\ell\rangle}}^{\tau_t^{\langle\ell\rangle}} \boldsymbol{E}^{\langle\ell\rangle}\mathrm{d}\xi_3 , \tag{39a}$$

$$\bar{\boldsymbol{B}} = \sum_{\ell=1}^{N_\ell} \int_{\tau_b^{\langle\ell\rangle}}^{\tau_t^{\langle\ell\rangle}} \boldsymbol{E}^{\langle\ell\rangle}\xi_3\mathrm{d}\xi_3 , \tag{39b}$$

$$\bar{\boldsymbol{D}} = \sum_{\ell=1}^{N_\ell} \int_{\tau_b^{\langle\ell\rangle}}^{\tau_t^{\langle\ell\rangle}} \boldsymbol{E}^{\langle\ell\rangle}\xi_3^2\mathrm{d}\xi_3 , \tag{39c}$$

$$\bar{\boldsymbol{C}} = \bar{\boldsymbol{B}}^T , \tag{39d}$$

where $\tau_b^{\langle\ell\rangle}$ and $\tau_t^{\langle\ell\rangle}$ are the values of $\xi_3$ at the bottom and top surfaces of the $\ell$-th layer, respectively, with the thickness of the $\ell$-th layer given by $\tau^{\langle\ell\rangle} = \tau_t^{\langle\ell\rangle} - \tau_b^{\langle\ell\rangle}$. These constitutive matrices in Voigt notation are used to obtain the constitutive tensors. Taking as example the first matrix $\bar{\boldsymbol{A}}$, the associated tensor is obtained through $\bar{\mathbb{A}}_{\alpha\beta\gamma\delta} = \bar{A}_{ab}$, using the correspondences $\alpha\beta \longleftrightarrow a$ and $\gamma\delta \longleftrightarrow b$, where the indices 11, 22, 12, and 21 correspond to 1, 2, 3, and 3, respectively. Finally, the laminate constitutive tensors in the local covariant basis are obtained trough the following transformation law:

$$\mathbb{A}^{\alpha_1\beta_1\gamma_1\delta_1} = \bar{\mathbb{A}}^{\alpha_2\beta_2\gamma_2\delta_2}\,(\boldsymbol{m}_{\alpha_2}\cdot\boldsymbol{a}^{\alpha_1})(\boldsymbol{m}_{\beta_2}\cdot\boldsymbol{a}^{\beta_1})(\boldsymbol{m}_{\gamma_2}\cdot\boldsymbol{a}^{\gamma_1})(\boldsymbol{m}_{\delta_2}\cdot\boldsymbol{a}^{\delta_1}) . \tag{40}$$

# C Kirchhoff-Love fluxes

The covariant and contravariant component of the vector $\boldsymbol{t}$ are obtained as $t_\alpha = \boldsymbol{t} \cdot \boldsymbol{a}_\alpha$ and $t^\alpha = \boldsymbol{t} \cdot \boldsymbol{a}^\alpha$, respectively. The same applies to the vector $\boldsymbol{n}$. The bending and twisting moments are derived from $M_{nn} = M^{\alpha\beta} n_\alpha n_\beta$ and $M_{nt} = M^{\alpha\beta} n_\alpha t_\beta$, respectively. The covariant derivative of the moment tensor is obtained as

$$M^{\alpha\beta}_{|\gamma} = M^{\alpha\beta}_{,\gamma} + \Gamma^\alpha_{\lambda\gamma} M^{\lambda\beta} + \Gamma^\beta_{\lambda\gamma} M^{\alpha\lambda} , \tag{41}$$

where $\Gamma^\gamma_{\alpha\beta}$ represents the Christoffel symbols of the second kind which are defined as

$$\Gamma^\gamma_{\alpha\beta} = \boldsymbol{a}^\gamma \cdot \boldsymbol{a}_{\alpha,\beta} . \tag{42}$$

Recalling Eq.(12), the coordinate derivative of the moment tensor is computed as

$$M^{\alpha\beta}_{,\rho} = \mathbb{C}^{\alpha\beta\gamma\delta}_{,\rho}\ \varepsilon_{\gamma\delta} + \mathbb{C}^{\alpha\beta\gamma\delta}\ \varepsilon_{\alpha\beta,\rho} + \mathbb{D}^{\alpha\beta\gamma\delta}_{,\rho}\ \kappa_{\gamma\delta} + \mathbb{D}^{\alpha\beta\gamma\delta}\ \kappa_{\alpha\beta,\rho} . \tag{43}$$

The derivatives appearing in Eq.(43) are not reported here for the sake of conciseness. However, their computation through the chain rule is straightforward, albeit somewhat laborious. To what regards the second term in Eq.(18b), the following relationship holds:

$$(M^{\alpha\beta} n_\alpha t_\beta)_{,t} = (M^{\alpha\beta})_{,t} n_\alpha t_\beta + M^{\alpha\beta} (n_\alpha)_{,t} t_\beta + M^{\alpha\beta} n_\alpha (t_\beta)_{,t} , \tag{44}$$

where the notation $(\bullet)_{,t}$ is used to denote the arc-length derivative along the curve that identifies $\boldsymbol{t}$. Let us suppose that this curve is known in the Euclidean space through a map $\boldsymbol{y} = \boldsymbol{y}(\tau)$, being $\tau$ an auxiliary curvilinear direction. Then, the tangent unit vector is computed as $\boldsymbol{t} = \frac{\boldsymbol{y}_{,\tau}}{|\boldsymbol{y}_{,\tau}|}$, while its derivative with respect to $\tau$ is

$$\boldsymbol{t}_{,\tau} = \frac{(\boldsymbol{I} - \boldsymbol{t} \otimes \boldsymbol{t})}{|\boldsymbol{y}_{,\tau}|} \boldsymbol{y}_{,\tau\tau} , \tag{45}$$

where $\boldsymbol{I}$ is the identity $3 \times 3$ tensor. The arc-length derivative of $\boldsymbol{t}$ is further obtained as

$$(\boldsymbol{t})_{,t} = \frac{\boldsymbol{t}_{,\tau}}{|\boldsymbol{y}_{,\tau}|} . \tag{46}$$

To what regards the normal vector $\boldsymbol{n} = \boldsymbol{t} \times \boldsymbol{a}_3$, its arc-length derivative is computed as

$$(\boldsymbol{n})_{,t} = (\boldsymbol{t})_{,t} \times \boldsymbol{a}_3 + \boldsymbol{t} \times (\boldsymbol{a}_3)_{,t} . \tag{47}$$

However, it is worth noting that when dealing with shells coupling, the relative orientation of $\boldsymbol{t}$, $\boldsymbol{n}$ and $\boldsymbol{a}_3$ may be such that $\boldsymbol{n} = \boldsymbol{a}_3 \times \boldsymbol{t}$, which implicates a straightfoward adaptation of Eq.(47). The arc-length derivative of $\boldsymbol{a}_3$ is computed as $(\boldsymbol{a}_3)_{,t} = \boldsymbol{a}_{3,\gamma} t^\gamma$, where, introducing the vector $\boldsymbol{p} = \boldsymbol{a}_1 \times \boldsymbol{a}_2$, $\boldsymbol{a}_{3,\gamma}$ is obtained as

$$\boldsymbol{a}_{3,\gamma} = \frac{(\boldsymbol{I} - \boldsymbol{a}_3 \otimes \boldsymbol{a}_3)}{\lambda} \boldsymbol{p}_{,\gamma} , \tag{48}$$

and $\boldsymbol{p}_{,\gamma}$ is computed through the chain rule. The arc-length derivatives of the covariant components of $\boldsymbol{t}$ and $\boldsymbol{n}$ are obtained as $(t_\alpha)_{,t} = (\boldsymbol{t})_{,t} \cdot \boldsymbol{a}_\alpha$ and $(n_\alpha)_{,t} = (\boldsymbol{n})_{,t} \cdot \boldsymbol{a}_\alpha$. Finally, the arc-length derivative of $M^{\alpha\beta}$ is obtained as $(M^{\alpha\beta})_{,t} = M^{\alpha\beta}_{|\gamma} t^\gamma$.